\documentclass{amsart}
\usepackage{latexsym}
\usepackage{amsmath,amsthm,amsfonts,amscd,eucal}
\usepackage{graphicx}
\usepackage{epsfig}

\numberwithin{equation}{section}

\def\ca{{\mathcal A}}

\def\cc{{\mathcal C}}

\def\cf{{\mathcal F}}
\def\cg{{\mathcal G}}

\def\cp{{\mathcal P}}

\def\car{{\mathcal R}}

\def\bc{{\mathbb C}}

\def\bn{{\mathbb N}}
\def\br{{\mathbb R}}
\def\bz{{\mathbb Z}}

\def\a{\alpha}

\def\g{\gamma}        \def\G{\Gamma}
\def\d{\delta}        \def\D{\Delta}
\def\eps{\varepsilon}

\def\th{\vartheta}    

\def\k{\kappa}
\def\l{\lambda}       \def\La{\Lambda}
\def\m{\mu}

\def\r{\rho}
\def\s{\sigma}        
     
\def\t{\tau}

        \def\O{\Omega}

\def\itm#1{\item{$(#1)$}}
\newcommand{\set}[1]{\left\{#1\right\}}
\newcommand{\card}[1]{\left| \set{#1} \right|}
\newcommand{\supp}{\text{supp}}
\newcommand{\conv}{\text{conv}\,}

\newcommand{\Det}{\text{det}}

\newcommand{\Ci}{\cc}  
\newcommand{\Nt}{\cc^{\mathrm{notail}}}  
\newcommand{\Ta}{\cc^{\mathrm{tail}}}  
\renewcommand{\Re}{\car}  
\renewcommand{\Pr}{\cp}  

\newtheorem{Thm}{Theorem}[section]
\newtheorem{Cor}[Thm]{Corollary}
\newtheorem{Prop}[Thm]{Proposition}
\newtheorem{Lemma}[Thm]{Lemma}

\theoremstyle{definition}
\newtheorem{Dfn}[Thm]{Definition}
\newtheorem{exmp}[Thm]{Example}

\theoremstyle{remark}
\newtheorem{rem}[Thm]{Remark}

\begin{document}

\title[A trace on fractal graphs and the Ihara zeta function]
{A trace on fractal graphs\\ and  the Ihara zeta function}%
\author{Daniele Guido, Tommaso Isola, Michel L. Lapidus}%
\address{(D.G., T.I.) Dipartimento di Matematica, Universit\`a di Roma
``Tor Vergata'', I--00133 Roma, Italy.}%
\email{guido@mat.uniroma2.it, isola@mat.uniroma2.it}%
\address{(M.L.L.) Department of Mathematics, University of California,
Riverside, CA 92521-0135, USA.}%
\email{lapidus@math.ucr.edu}%
\thanks{The first and second authors were partially
supported by MIUR, GNAMPA and by
the European Network ``Quantum Spaces - Noncommutative Geometry"
HPRN-CT-2002-00280. The third author was partially supported
by the National Science Foundation, the Academic Senate of the 
University of California, and GNAMPA}%

\begin{abstract}
    Starting with Ihara's work in 1968, there has been a growing
    interest in the study of zeta functions of finite graphs, by
    Sunada, Hashimoto, Bass, Stark and Terras, Mizuno and Sato, to
    name just a few authors.  Then, Clair and Mokhtari-Sharghi have
    studied zeta functions for infinite graphs acted upon by a
    discrete group of automorphisms.  The main formula in all these
    treatments establishes a connection between the zeta function,
    originally defined as an infinite product, and the Laplacian of
    the graph.  In this article, we consider a different class of
    infinite graphs.  They are fractal graphs, i.e. they enjoy a
    self-similarity property.  We define a zeta function for these
    graphs and, using the machinery of operator algebras, we prove a
    determinant formula, which relates the zeta function with the
    Laplacian of the graph.  We also prove functional equations, and a
    formula which allows approximation of the zeta function by the
    zeta functions of finite subgraphs.
\end{abstract}

\maketitle

\setcounter{section}{-1}

\section{Introduction}

The Ihara zeta function, originally associated to certain groups and
then combinatorially reinterpreted as associated with finite graphs or
with their infinite coverings, is defined here for a new class of
infinite graphs, called self-similar graphs.  The corresponding
determinant formula and functional equations are established.

The combinatorial nature of the Ihara zeta function was first observed
by Serre (see \cite{Serre}, Introduction), but it was only through the
works of Sunada \cite{Sunada}, Hashimoto \cite{HaHo,Hashi} and Bass
\cite{Bass} that it became a graph-theoretical object, at the same
time keeping some number-theoretically flavoured properties, like the
Euler product formula or the functional equation.

The Ihara zeta function \cite{Ihara} was written as an infinite
product (Euler product) over $G$-conjugacy classes of primitive
elements in a group $G$, namely elements whose centralizer (in $G$) is
generated by the element itself.  As explained in detail in the
introductions of \cite{Bass,StTe}, Ihara's construction can be
rephrased in terms of a regular (i.e. constant number of edges
spreading from each vertex) finite graph $X$, its universal covering
$Y$ and the corresponding structure group $G=\pi_{1}(X)$.  By the
homotopic nature of $G$, one may equivalently represent $G$-conjugacy
classes in terms of suitably reduced primitive cycles on the graph
$X$.  Here, a reduced cycle on $X$ (of length $m$) is a set
$\{e_{j},j\in\bz_{m}\}$ where the starting vertex of $e_{j+1}$
coincides with the ending vertex of $e_j$, and $e_{j+1}$ is not the
opposite of $e_j$, for $j\in\bz_{m}$.  Besides, a cycle is not
primitive if it is obtained by repeating the same cycle more than
once.  Finally, denoting by $\Pr$ the set of (reduced) primitive
cycles and by $|C|$ the length of a cycle $C$, the Ihara zeta function
of a finite graph can be written, for $|u|$ small enough, as
\begin{equation}\label{finiteIhara}
    Z(u) := \prod_{C\in \Pr}(1-u^{|C|})^{-1}.
\end{equation}

The Ihara zeta function may be considered as a modification of the
Selberg zeta function, cf.  e.g. \cite{Bass}, and was originally
written in terms of the variable $s$, the relation with $u$ being, for
$(q+1)$-regular graphs, $u=q^{-s}$: $Z(s) = \prod_{C\in
\Pr}(1-(q^{|C|})^{-s})^{-1}$.  In this form, the relation with the
Riemann zeta function is apparent.  The latter, one of the primary
examples of number-theoretic zeta functions, may indeed be expanded as
an Euler product $\zeta(s)=\prod_{p}(1-p^{-s})^{-1}$ for $Re\, s>1$,
where $p$ ranges over all rational primes.  Further, $\zeta(s)$ may be
meromorphically extended to the whole complex plane, where its
completion $\xi(s)=\pi^{-s/2}\G(s/2)\zeta(s)$ verifies the functional
equation $\xi(s)=\xi(1-s)$.

Concerning the Ihara zeta function, many efforts have been exerted to
clarify one of its main properties, namely the fact that its inverse
is a polynomial, more precisely the determinant of a matrix-valued
polynomial.  Such equality, known as the determinant formula, has been
proved by Bass \cite{Bass} for non-regular graphs, and reproved by
several authors \cite{StTe,FoZe,Bartholdi,No2,KoSu}, different proofs
corresponding to some generalizations of the original setting and also
to a shift in the techniques, from purely algebraic to more
combinatorial and functional analytic.  As a consequence of the
determinant formula, the zeta function of Ihara meromorphically
extends to the whole complex plane, and its completions satisfy a
functional equation.  Let us mention that other number-theoretic
properties, like a version of the Riemann hypothesis, have been
studied for the Ihara zeta function, the graphs satisfying it being
completely characterized, see \cite{StTe} for a simple proof.  More
applications of the Ihara zeta function are contained in
\cite{Bass,Hashi1,Hashi2,Serre1,No1,StTe2,StTe3,HSTe}.  Different
generalizations of the Ihara zeta function are considered in the
literature, see \cite{Bartholdi,MiSa} and the references therein.

As for index theorems and geometric invariants, where the theory was
extended from compact manifolds to covering manifolds by Atiyah
\cite{Atiyah}, it was observed by Clair and Mokhtari-Sharghi
\cite{ClMS1} that Ihara's construction can be extended to infinite
graphs on which a group $\G$ acts isomorphically and with finite
quotient.

Such a group $\G$ gives rise to an equivalence relation between
(primitive) cycles: to obtain the Ihara zeta function for coverings,
formula (\ref{finiteIhara}) should be modified in the sense that the
Euler product is taken over $\G$-equivalence classes of primitive
cycles, and each factor should be normalized by an exponent related
with the cardinality of the stabilizer of the cycle.  This subject has
been implicitely considered by Bass \cite{Bass} and then extensively
studied by Clair and Mokhtari-Sharghi in \cite{ClMS1,ClMS2}.  In
particular, a determinant formula has been established, this result
being deduced as a specialization of the treatment of group actions on
trees (the so-called theory of tree lattices, as developed by Bass,
Lubotzky and others, see \cite{BaLu}).

Moreover, as for the index theorem for coverings, the von Neumann
$\G$-trace on periodic operators acting on the graph plays an
important role.  A self-contained proof of the main results for the
Ihara zeta function of covering graphs, together with the proof of a
conjecture of Grigorchuk and $\dot{\text{Z}}$uk \cite{GrZu}, is
contained in \cite{GILa03}; cf.  also \cite{GILa02} for an
introduction to the subject.

Continuing the parallel with index theorems and geometric invariants,
a step beyond the periodic case was taken by Roe \cite{Roe}, who
proved an index theorem for open manifolds, and by Farber
\cite{Farber}, who proposed a general approximation scheme for
$L^{2}$-invariants.  In the same spirit, Grigorchuk and
$\dot{\text{Z}}$uk \cite{GrZu} proposed an approximation scheme for
Ihara zeta functions of infinite graphs.

In this paper we intend to realize this proposal for a suitable family
of infinite graphs, the self-similar graphs.  This family contains
many examples of what are known as fractal graphs in the literature,
see \cite{Barl,HaKu}.  These include the Gasket, Vicsek, Lindstrom and
Carpet graphs, see figures \ref{fig:GasketVicsek},
\ref{fig:LindstromCarpet}.
 
On the one hand, these graphs can be appoximated by finite graphs as
in the case of amenable coverings, namely the ratio between the size
of the boundary and the total size of the approximating finite graphs
becomes smaller and smaller.  On the other hand, they possess many
local isomorphisms, whose domains correspond to arbitrarily large
portions of the graph.  These local isomorphisms guarantee that the
approximation works, as we explain below.  We define the Ihara zeta
function of a self-similar graph as an Euler product over equivalence
classes of primitive cycles, the equivalence being given by local
isomorphisms, and each factor having a normalization exponent related
with the average multiplicity of the given equivalence class, cf. 
Definition \ref{Dfn:Zeta}.  The existence of such an average
multiplicity is the first consequence of self-similarity.

In order to prove a determinant formula, we first define a
C$^{*}$-algebra of operators acting on $\ell^{2}$ of the vertices of
the graph, containing in particular the adjacency operator and the
degree operator, and then define a normalized trace on its elements,
given by the limit of the normalized traces of the restriction of an
operator to the approximating finite graphs.  The existence of such a
trace is another important consequence of self-similarity.

Let us notice that the existence of a trace for finite propagation
operators on spaces with an amenable exhaustion was first established in
\cite{Roe}; see also \cite{GuIs07}. However, these traces depended on a
generalized limit procedure.  The independence of the choice of such a 
generalized limit was proved in \cite{CGIs01} for self-similar CW-complexes 
and in \cite{Elek1} for abstract quasicrystal graphs.

We then define an analytic determinant for C$^*$-algebras with a
finite trace.  In contrast with the Fuglede--Kadison determinant
\cite{FuKa} for finite von Neumann algebras, our determinant depends
analytically on its argument, but, as a drawback, is defined on a
smaller domain and obeys weaker properties.  Such an analytic
determinant, based on the trace described above, is used in the
determinant formula.  As for the finite graph case, the determinant
formula allows the extension of the zeta function to a larger domain,
and finally implies, in the case of regular graphs, the validity of
the functional equation for suitable ``completions'' of the zeta
function.  As a final check for the approximation structure, we show
that the zeta functions of the approximating finite graphs, with a
suitable renormalizing exponent, converge to the zeta function of the
self-similar graph.

The structure of the paper is the following.  After having introduced
some preliminary notions, we define in Section \ref{sec:selfSimilar}
the class of graphs we are interested in, and, in the following two
sections, we construct a trace on the C$^{*}$-algebra of geometric
operators on the graph, and a determinant, on a suitable subclass of
operators, extending previous work of the authors \cite{GILa02}.

Then, after some technical preliminaries in Section
\ref{sec:PrimeCycles}, we define, in Section \ref{sec:Zeta}, our zeta
function, and show that it is a holomorphic function.  In Section
\ref{sec:detFormula}, we prove a corresponding determinant formula. 
In Section \ref{sec:essRegGraph}, we establish several functional
equations, for different completions of the zeta function.
 
Finally, in Section \ref{sec:approx}, we prove that the zeta function
is the limit of a sequence of (appropriately normalized) zeta
functions of finite subgraphs, which shows that our definition of the
zeta function is a natural one.
 
In closing this introduction, we note that, in \cite{LavF1,LavF2}, a
variety of zeta functions are associated to certain classes of
(continuous rather than discrete) fractals.  They are, however, of a
very different nature than the Ihara-type zeta functions of fractal
graphs considered in this paper.
  
The contents of this paper have been presented at the 2006 Spring
Western Section Meeting of the American Mathematical Society in San
Francisco in April 2006, and at the $21^{st}$ conference on Operator
Theory in Timisoara (Romania) in July 2006.

\section{Preliminaries}

 In this section, we recall some terminology from graph theory, and 
 introduce the class of geometric operators on an infinite graph. 

 A {\it simple graph} $X=(VX,EX)$ is a collection $VX$ of objects,
 called {\it vertices}, and a collection $EX$ of unordered pairs of
 distinct vertices, called {\it edges}.  The edge $e=\set{u,v}$ is
 said to join the vertices $u,v$, while $u$ and $v$ are said to be
 {\it adjacent}, which is denoted $u\sim v$.  A {\it path} (of length
 $m$) in $X$ from $v_0\in VX$ to $v_m\in VX$, is
 $(v_{0},\ldots,v_{m})$, where $v_{i}\in VX$, $v_{i+1}\sim v_i$, for
 $i=0,...,m-1$ (note that $m$ is the number of edges in the path).  A
 path is {\it closed} if $v_{m}=v_{0}$.
 
 We assume that $X$ is countable and connected, $i.e.$ there is a path
 between any pair of distinct vertices.  Denote by $deg(v)$ the degree
 of $v\in VX$, $i.e.$ the number of vertices adjacent to $v$.  We
 assume that $X$ has bounded degree, $i.e.$ $d := \sup_{v\in VX}
 deg(v) <\infty$.  Denote by $\r$ the combinatorial distance on $VX$,
 that is, for $v,w\in VX$, $\r(v,w)$ is the length of the shortest
 path between $v$ and $w$.  If $\Omega \subset VX$, $r\in\bn$, we
 write $B_r(\O) := \cup_{v\in\O} B_r(v)$, where $B_{r}(v) := \set{
 v'\in VX: \r(v',v)\leq r}$.
 
 \begin{Dfn}[Finite propagation operators]
    A bounded linear operator $A$ on $\ell^2 (VX)$ has {\it finite
    propagation} $r=r(A)\geq 0$ if, for all $v\in VX$, we have $\supp
    (Av)\subset B_{r}(v)$ and $\supp (A^{*}v)\subset B_{r}(v)$, where
    $A^*$ is the Hilbert space adjoint of $A$.
 \end{Dfn}
 
 \begin{rem}
     \itm{i} Finite propagation operators have been called ``bounded
     range" by other authors.  
     
     \itm{ii} Let us note that, if $\lambda \in \bc$ and $A, B$ are
     finite propagation operators,
     $$
     r(\lambda A + B)= r(A) \vee r(B)\, ,\qquad r(AB) = r(A) + r(B),
     $$
     showing that finite propagation operators form a $^*$-algebra.
 \end{rem}

 \begin{Dfn}[Local Isomorphisms and Geometric Operators]
     A {\it local isomorphism} of the graph  $X$ is a triple
     \begin{equation}
	 \Bigl(S(\gamma)\, ,R(\gamma)\, ,\gamma \Bigr),
     \end{equation}
     where $S(\gamma)\, ,R(\gamma)$ are subgraphs of $X$ and
     $\gamma : S(\gamma)\to R(\gamma)$ is a graph isomorphism. 

     \noindent The local isomorphism $\gamma$ defines a {\it partial
     isometry} $U(\gamma) : \ell^2 (VX)\to \ell^2 (VX)$, by setting
     \begin{align*}
	 U(\g)(v):=
	 \begin{cases}
	     \g(v)& v\in V(S(\g)) \\
	     0& v\not\in V(S(\g)), 
	 \end{cases}
     \end{align*}
     and extending by linearity.  A bounded operator $T$ acting on
     $\ell^2 (VX)$ is called {\it geometric} if there exists $r\in\bn$
     such that $T$ has finite propagation $r$ and, for any local
     isomorphism $\gamma$, any $v\in VX$ such that $B_{r}(v)\subset
     S(\g)$ and $B_{r}(\g v)\subset R(\g)$, one has
     \begin{equation}\label{eq:geometric}
       TU(\gamma)v = U(\gamma)Tv,\quad T^*U(\gamma)v = U(\gamma)T^*v\, .
     \end{equation}     
 \end{Dfn}
 
 \begin{rem}\label{rem:degree}
     A local isomorphism $\g$ does not necessarily preserve the degree
     of a vertex.  It does, however, preserve the degree of any vertex
     $v\in S(\g)$ such that $B_{1}(v)\subset S(\g)$ and $B_{1}(\g
     v)\subset R(\g)$.
 \end{rem}

 Recall that the {\it adjacency matrix} of $X$,
 $A=\big(A(v,w)\big)_{v,w\in VX}$, and the {\it degree matrix} of $X$,
 $D=\big(D(v,w)\big)_{v,w\in VX}$ are defined by
 \begin{equation}\label{eq:adjacency}
     A(v,w)=
     \begin{cases} 
	 1&v\sim w\\
	 0&\text{otherwise} 
     \end{cases}
 \end{equation}
 and
 \begin{equation}\label{eq:degree}
     D(v,w)=
     \begin{cases} 
	 \deg(v)&v= w\\
	 0&\text{otherwise.} 
     \end{cases}
 \end{equation}
 
 \begin{Prop}\label{Prop:3.7}
     Geometric operators form a $^*$-algebra containing the adjacency
     operator $A$ and the degree operator $D$.
 \end{Prop}
 \begin{proof}
     The set of geometric operators is clearly a vector space which is
     $^*$-closed.  Let us now prove that it is also closed with
     respect to the product.  Let $T_1$ and $T_2$ be geometric
     operators, let $\g$ be a local isomorphism and let $r_1>0$ (resp. 
     $r_2>0$) be such that (\ref{eq:geometric}) holds for $T_1$ (resp. 
     $T_2$).  Let $r:=r_1+r_2$.  Then, for any $v\in VX$ for which
     $B_{r}(v)\subset S(\g)$ and $B_{r}(\g v)\subset R(\g)$, one has
     $$
     [T_1T_2,U(\g)]v = T_1[T_2,U(\g)]v + [T_1,U(\g)] T_2v = 0,
     $$
     where $[\ \cdot\ ,\ \cdot\ ]$ deonotes the commutator, because
     $[T_2,U(\g)]v=0$ and $T_2v = \sum_{j=1}^k c_j w_j$ is a linear
     combination of vertices $w_j$ belonging to $B_{r_2}(v)$, so that
     $B_{r_1}(w_j)\subset B_r(v)\subset S(\g)$ and $B_{r_1}(\g
     w_j)\subset B_r(\g v)\subset R(\g)$ for all $j=1,...,k$.  Hence
     $[T_1,U(\g)] T_2v = 0$ by linearity.  An analogous argument shows
     that $[(T_1T_2)^*,U(\g)]v=0$.
     
     We now prove that $A$ is geometric.  First recall that $\|A\|\leq
     d <\infty$ (see \cite{Mohar}, and {\it e.g.} \cite{MoWo}).  Let
     $\g$ be a local isomorphism and let $v\, ,v^{\prime}\in VX$ be
     such that $B_{1}(v)\subset S(\g)$ and $B_{1}(\g v)\subset R(\g)$. 
     Then, if $v'\not\in R(\g)$, because $\supp(A v)\subset B_{1}(v)
     \subset S(\g)$ and $\supp(A(\g v))\subset B_{1}(\g v) \subset
     R(\g)$, we obtain
     $$
     (AU(\g)v, v') = 0 = (U(\g)Av, v').
     $$
     Thus, let us suppose that $v'\in R(\g)$, so that $v'=\g v''$, for
     $v''\in S(\g)$.  Then $v' = \g v'' \sim \g v \iff v'' \sim v$,
     because $\g$ is a local isomorphism, so that
     \begin{align*}    
	 (AU(\g)v , v' ) & =
	 (A(\g v) , \g v'' )  = (Av,v'') \\
	 & = (Av, U(\g)^{*}v')  = (U(\g)Av, v').
     \end{align*}
     By linearity, we deduce that $A$ is geometric.
     
     Finally, the degree operator $D$ is a multiplication operator,
     hence it has zero propagation.  However, by Remark
     \ref{rem:degree}, it is geometric with constant $r=1$.
 \end{proof}
 
\section{Self-similar graphs}\label{sec:selfSimilar}

 In this section, we introduce the class of self-similar graphs.  This
 class contains many examples of what are usually called fractal
 graphs, see $e.g.$ \cite{Barl,HaKu}.

 If $K$ is a subgraph of $X$, we call {\it frontier} of $K$, and
 denote by $\cf(K)$, the family of vertices in $V K$ having distance 1
 from the complement of $V K$ in $VX$.

 \begin{Dfn}[Amenable graphs]
     A countably infinite graph with bounded degree $X$ is {\it
     amenable} if it has an {\it amenable exhaustion}, namely, an
     increasing family of finite subgraphs $\{K_n : n\in \mathbb{N}\}$
     such that $\cup_{n\in\bn} K_n = X$ and
    \begin{equation*}
	\frac{|\cf(K_n)|}{|K_n|}\to 0\qquad {\rm as}\,\,\, n\to
	\infty\, ,
    \end{equation*}
    where $|K_{n}|$ stands for $|VK_{n}|$ and $|\cdot|$ denotes the 
    cardinality.
\end{Dfn}

\begin{Dfn}[Self-similar graphs] \label{def:Quasiperiodic} 
    A countably infinite graph with bounded degree $X$ is
    {\it self-similar} if it has an amenable exhaustion $\{K_{n}\}$
    such that the following  conditions  $(i)$ and $(ii)$ hold: 
    
    \itm{i} For every $n\in\bn$, there is a finite set of local
    isomorphisms $\cg(n,n+1)$ such that, for all $\gamma\in
    \cg(n,n+1)$, one has $S(\gamma) = K_n$,
    \begin{equation}
	\bigcup_{\gamma\in \cg(n,n+1)} \gamma(K_n) = K_{n+1},
    \end{equation}
    and moreover, if $\gamma , \gamma^\prime\in \cg(n,n+1)$ with
    $\gamma\neq\gamma^\prime$,
    \begin{equation}
	V( \gamma K_n ) \cap V( \gamma' K_n ) = 
	\cf(\gamma K_n ) \cap \cf( \gamma' K_n ).
    \end{equation}
    
    \itm{ii} We then define $\cg(n,m)$, for $n<m$, as the set of all
    admissible products $\g_{m-1}\cdot\dots\cdot\g_{n}$, $\g_{i}\in
    \cg(i,i+1)$, where ``admissible'' means that, for each term of the
    product, the range of $\g_{j}$ is contained in the source of
    $\g_{j+1}$.  We also let $\cg(n,n)$ consist of the identity
    isomorphism on $K_{n}$, and $\cg(n):=\cup_{m\geq n}\cg(n,m)$.  We
    can now define the $\cg$-{\it invariant frontier} of $K_{n}$:
    $$
    \cf_{\cg}(K_{n})=
    \bigcup_{\gamma\in \cg(n)}\g^{-1}\cf( \gamma K_{n} ),
    $$
    and we require that
    \begin{equation}\label{strongreg}
	\frac{|\cf_{\cg}(K_n)|}{|K_n|}\to 0\qquad {\rm
	as}\; n\to \infty\, .
    \end{equation}
\end{Dfn}

 In the rest of the paper, we denote by $\cg$ the family of all
 local isomorphisms which can be written as (admissible) products
 $\g_1^{\eps_1} \g_2^{\eps_2}...\g_k^{\eps_k}$, where
 $\g_i\in\cup_{n\in\bn} \cg(n)$, $\eps_i\in\set{-1,1}$, for
 $i=1,...,k$ and $k\in\bn$.
 
 \begin{rem}
     $(i)$ Condition $(i)$ of the Definition above means that each
     $K_{n+1}$ is given by a finite union of copies of $K_n$, and
     these copies can only have the frontier in common.  In particular
     the number of such copies may vary with $n$, as well as the
     intersections of the frontiers.  One may also read this as a
     constructive recipe: choose a finite graph $K_1$, then arrange
     finitely many copies of it by assigning possible intersections
     and call this new graph $K_2$.  Now repeat the operation with
     $K_2$, and so on.  No finite requirement is needed; one should
     only guarantee that the degree remains bounded and that condition
     (\ref{strongreg}) is satisfied.
 
     \itm{ii} In the examples described below the map from $K_n$ to
     $K_{n+1}$ is essentially the same for all $n$, and is related to
     the construction of a self-similar fractal.  One may generalize
     this construction to translation fractals, cf.  \cite{GuIs10}.
 
     \itm{iii} Another axiomatic construction of a family of
     self-similar graphs was considered by Kr\"{o}n and Teufl
     \cite{Kr0,Kr,KT}.  Their idea is based on the existence of a sort
     of dilation map on the vertices of the graph, and a procedure
     associating to the ``dilated vertices" a graph structure which
     reproduces exactly the original graph.  Their family is smaller
     than ours, but allowed them to make effective computations
     concerning the random walk and the asymptotic dimension of the
     graph.
 \end{rem}
 
 \begin{exmp}\label{ex:Graphs}
     Several examples of self-similar graphs are shown in figures
     \ref{fig:GasketVicsek}a, \ref{fig:GasketVicsek}b, 
     \ref{fig:LindstromCarpet}a, \ref{fig:LindstromCarpet}b.  They are
     called the Gasket graph, the Vicsek graph, the Lindstrom graph
     and the Carpet graph, respectively.
     \begin{figure}[ht]
 	 \centering
	 \psfig{file=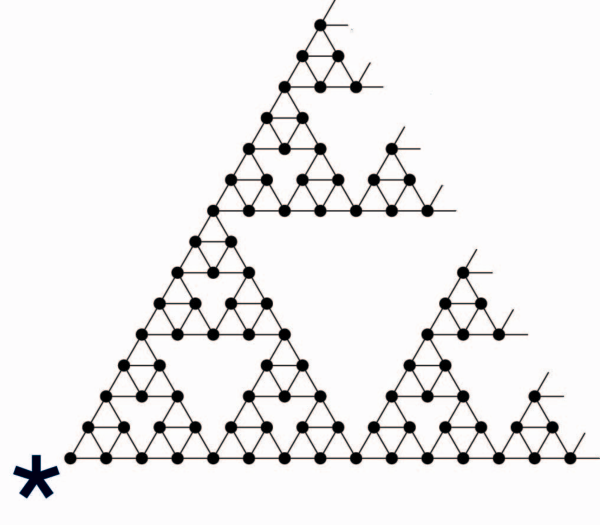,height=1.5in} \qquad \qquad
	 \psfig{file=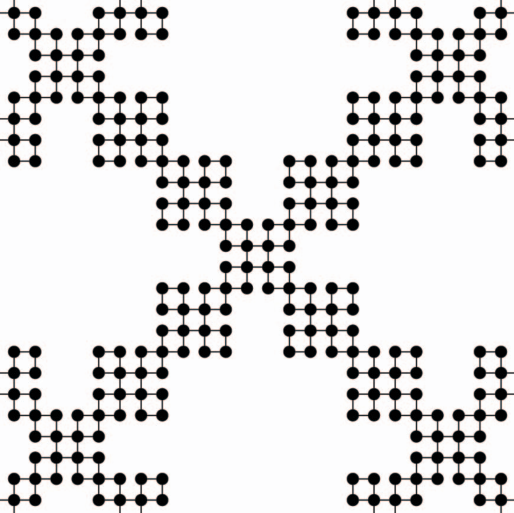,height=1.5in}
	 \caption{(a) Gasket graph. \qquad \qquad (b) Vicsek graph.}
	 \label{fig:GasketVicsek}
     \end{figure}
     \begin{figure}[ht]
	 \centering
	 \psfig{file=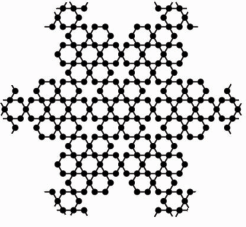,height=1.5in} \qquad \qquad
	 \psfig{file=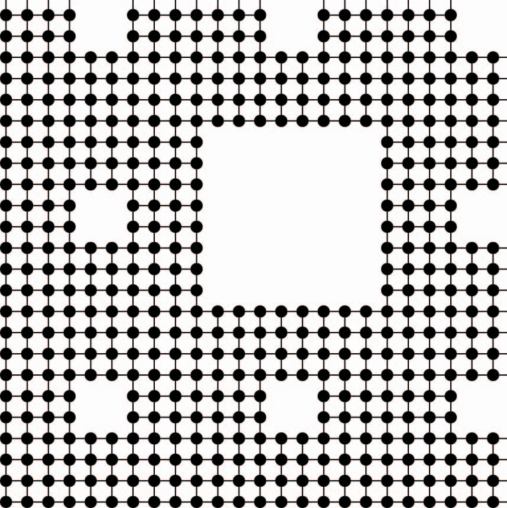,height=1.5in}
	 \caption{(a) Lindstrom graph. \qquad \qquad (b)  Carpet graph.}
	 \label{fig:LindstromCarpet}
     \end{figure}
 \end{exmp}
 
 These examples can all be obtained from the following general
 procedure (for more details, see \cite{CGIs01}, where the more
 general case of self-similar CW-complexes has been studied).  Assume
 that we are given a self-similar fractal in $\br^{p}$ determined by
 similarities $w_{1},\ldots,w_{q}$, with the same similarity
 parameter, and satisfying the Open Set Condition for a bounded open
 set whose closure is a convex $p$-dimensional polyhedron $\cp$, and
 let $K_1$ be the graph consisting of the vertices and edges of $\cp$. 
 If $\sigma = (\sigma_{1},..., \sigma_{n})$ is a multi-index of length
 $n$, we set $w_{\s}:= w_{\s_{n}}\circ\cdots\circ w_{\s_{1}}$, and
 assume that $w_{\s}\cp\cap w_{\s'}\cp$ is a (facial) subpolyhedron of
 both $w_{\s}\cp$ and $w_{\s'}\cp$, with $|\s|=|\s'|$.  Finally, we
 choose an infinite multi-index $\t$ and set $K_{n+1}:= \cup_{|\s'|=n}
 w_{\t|_{n}}^{-1}w_{\s'} K_1$, where $\t|_n$ is the multi-index of
 length $n$ obtained by truncation of $\t$ to its first $n$ letters. 
 Then $X := \cup_{n=1}^\infty K_n$ is a self-similar graph, with
 amenable exhaustion $\set{K_n}$ and family of local isomorphisms
 given by $\cg(n,n+1) := \set{\g_{i} :=
 w_{\t|_{n+1}}^{-1}w_{i}w_{\t|_n}: i=1,...,q}$, $n\in \bn$.  Observe
 that, in the above example \ref{ex:Graphs}, the generating polyhedron
 is an equilateral triangle, in the case of the Gasket graph, a
 square, in the case of the Vicsek or the Carpet graphs, and a regular
 hexagon, in the case of the Lindstrom graph.

\section{A trace on geometric operators}

 In this section, we construct a trace on the algebra of geometric
 operators on a self-similar graph.  In the following section, this
 trace will be used to define a determinant on some class of operators
 on the graph.

\begin{Thm}\label{thm:trace}
    Let $X$ be a self-similar graph, and let $\ca (X)$ be the
    C$^*$-algebra defined as the norm closure of the $^*$-algebra of
    geometric operators.  Then, on $\ca(X)$, there is a well-defined
    trace state $Tr_{\cg}$ given by
    \begin{equation}
	Tr_{\cg} (T) = \lim_n \frac{Tr\bigl( P(K_n) T \bigr)}
	{Tr\bigl( P(K_n) \bigr)},
    \end{equation}
    where $P(K_n)$ is the orthogonal projection of $\ell^{2}(VX)$ onto
    its closed subspace $\ell^{2}(V K_n)$.
\end{Thm}
\begin{proof}
    For a finite subset $N\subset VX$, denote by $P(N)$ the orthogonal 
    projection of $\ell^{2}(VX)$ onto $\text{span}(N)$.  Let
    us observe that, since $N$ is an orthonormal basis for $\ell^2
    (N)$, we have $Tr \bigl(P(N) \bigr) = |N|$.  For brevity's sake,
    we also use the notation $|K_{n}|:= |VK_{n}|$.
    
    \medskip
    \noindent {\bf First step}: some combinatorial results.
    
    a) Recall that $\displaystyle d=\sup_{v\in VX}| \set{v'\in VX:
    v'\sim v}|$, so that $\displaystyle \sup_{v\in VX}|B_1(v)|= d+1$.
    
    Then, since 
    $$
    B_{r+1}(v)=\bigcup_{v'\in B_{r}(v)} B_1(v),
    $$
    we get $|B_{r+1}(v)|\leq (d+1) |B_{r}(v)|$, giving
    $|B_{r}(v)|\leq (d+1)^{r}$, $\forall v\in VX$, $r\geq 0$.
    As a consequence, for any finite set $\Omega\subset VX$, 
    we have 
    $B_{r}(\Omega)=\displaystyle \bigcup_{v'\in \Omega}B_{r}(v')$,
    giving
    \begin{equation}\label{subsets}
	|B_{r}(\Omega)|\leq |\Omega| (d+1)^{r},\quad\forall r\geq 0.
    \end{equation}
    
    b) Let us set $\O_{n,r}=VK_{n}\setminus 
    B_{r}(\cf_{\cg}(K_{n}))$. Then, for any $\g\in\cg(n)$, we have
    \begin{equation}
	\g\O_{n,r}\subset \g VK_{n} = V(\g K_{n})\subset\g\O_{n,r}
	\cup B_{r}(\cf_{\cg}(\g K_{n})).
    \end{equation}
    Now assume $r\geq1$. Then the $\g\O_{n,r}$'s are disjoint for 
    different $\g$'s in $\cg(n,m)$. Therefore
    \begin{align}
	|K_{n}| &\leq |\O_{n,r}| + |\cf_{\cg}(K_{n})| (d+1)^{r},
	\label{ineqA}\\
	\left|VK_{m}\setminus\bigcup_{\g\in\cg(n,m)}\g\O_{n,r}\right| 
	&\leq |\cg(n,m)|\,|\cf_{\cg}(K_{n})| (d+1)^{r},\label{ineqB}\\
	|\cg(n,m)|\,|\O_{n,r}| & \leq |K_{m}|
	\leq |\cg(n,m)|\,|K_{n}|.\label{ineqC}
    \end{align}
    Indeed, (\ref{ineqA}) and (\ref{ineqC}) are easily verified, while
    \begin{align*}	    
	\left|VK_{m}\setminus\bigcup_{\g\in\cg(n,m)}\g\O_{n,r}\right|
	& = \left|\bigcup_{\g\in\cg(n,m)}\g
	VK_{n}\setminus\bigcup_{\g\in\cg(n,m)}\g\O_{n,r}\right|
	\\
	& \leq \sum_{\g\in\cg(n,m)} \left|\g 
	[VK_{n}\setminus \O_{n,r}]\right| \\
	& \leq |\cg(n,m)| \left| B_{r}(\cf_{\cg}(K_{n}))\right| \\
	& \leq |\cg(n,m)| \left| \cf_{\cg}(K_{n})\right| (d+1)^{r}.
    \end{align*}
     
    c) Set $\eps_{n}=\displaystyle{\frac{|\cf_{\cg}(K_{n})|}
    {|K_{n}|}}$ and recall that, by assumption, $\eps_{n}\to0$. 
    Putting together (\ref{ineqA}) and (\ref{ineqC}), we get
    $$
    |\cg(n,m)|\,|K_{n}| - |\cg(n,m)|\,|\cf_{\cg}(K_{n})| (d+1)^{r}
    \leq |K_{m}| \leq |\cg(n,m)|\,|K_{n}|,
    $$
    which implies
    \begin{equation}
	1-\eps_{n}(d+1)^{r}\leq \frac{|K_{m}|} {|\cg(n,m)|\, |K_{n}|}
	\leq 1.\label{ineqD}
    \end{equation}
    Choosing $n_{0}$ such that for all $n>n_{0}$,
    $\eps_{n}(d+1)^{r}\leq1/2$, we obtain
    \begin{equation}\label{ineqE}
	0\leq \frac{|\cg(n,m)|\, |K_{n}|}{|K_{m}|} -1 \leq 
	2\eps_{n}(d+1)^{r}\leq 1.
    \end{equation}
    Therefore, we deduce from (\ref{ineqB}) that
    \begin{align}\label{ineqF}
	\left|VK_{m}\setminus\bigcup_{\g\in\cg(n,m)}\g\O_{n,r}\right|
	& \leq |\cg(n,m)| \left| \cf_{\cg}(K_{n})\right| (d+1)^{r} 
	\\
	& = |\cg(n,m)| \left| K_{n}\right| \eps_{n}(d+1)^{r} \leq
	2 \left| K_{m}\right| \eps_{n}(d+1)^{r}. \notag
    \end{align}
    \\
    {\bf Second step}: the existence of the limit for geometric 
    operators.
    
    a) By definition of $U(\gamma)$, we have, for $\gamma\in
    \cg(n,m)$, with $n<m$,
    \begin{equation}
	U^*(\g) U(\g) = P(VK_n), \qquad U(\g) U^*(\g) = P(\g VK_n).
    \end{equation}
    Assume now that $T$ is a geometric
    operator on $\ell^2 (VX)$with finite propagation $r$.
    Then $T U(\g) P(\O_{n,r})= U(\g) T P(\O_{n,r})$ and
    $P(\g \O_{n,r}) = U(\g) P(\O_{n,r}) U(\g)^*$.  Hence,
	\begin{align}\label{eq:5.6}
	    Tr\bigl(T P(\g\O_{n,r}) \bigr) 
	    &= Tr\bigl(T U(\gamma) P(\O_{n,r}) U(\gamma)^* \bigr)\\
	    &= Tr\bigl(U(\gamma) T P(\O_{n,r}) U(\gamma)^* \bigr)\notag\\
	    &= Tr \bigl(TP(\O_{n,r}) U(\gamma)^* U(\gamma) \bigr) \notag\\
	    &= Tr \bigl(TP(\O_{n,r}) P(VK_n) \bigr) = Tr
	    \bigl(TP(\O_{n,r}) \bigr) .\notag
	\end{align}
    
    b) Let us show that the sequence is Cauchy:
    \begin{align*}
	&\left|\frac{Tr TP(VK_{n})}{Tr P(VK_{n})} -
	\frac{Tr TP(VK_{m})}{Tr P(VK_{m})} \right|\leq \\
	 &\leq \frac{|Tr T(P(VK_{n}) - P(\O_{n,r}))|} {|K_{n}|} + \frac
	{|Tr T(P(VK_{m}) - P(\cup_{\g\in\cg(n,m)}\g\O_{n,r}))|}
	{|K_{m}|} \\
	&\qquad+ \left|\frac{Tr T P(\O_{n,r})}{|K_{n}|} -
	\frac{|\cg(n,m)|\, |K_{n}|}{|K_{m}|} \frac{Tr T
	P(\O_{n,r})}{|K_{n}|} \right|\\
	 &\leq \|T\|\left( \frac{|VK_{n} \setminus \O_{n,r}|}{|K_{n}|}
	+ \frac{|VK_{m} \setminus \cup_{\g\in\cg(n,m)}\g\O_{n,r}|}
	{|K_{m}|} + \left|1 - \frac{|\cg(n,m)|\, |K_{n}|}{|K_{m}|}
	\right|\right)\\
	& \leq 5\|T\|\eps_{n}(d+1)^{r},
    \end{align*}
    where we used (\ref{eq:5.6}), in the first inequality, and
    (\ref{subsets}), (\ref{ineqF}), (\ref{ineqE}), in the third
    inequality. 
    \medskip
    
    \noindent {\bf Third step}: $Tr_{\cg}$ is a state on $\ca(X)$. 
    
    a) Let $T\in\ca(X)$, $\eps>0$.  Since $\ca(X)$ is the norm closure
    of the $^*$-algebra of geometric operators, we can find a
    geometric operator $T'\in\ca_{g}(X)$ such that
    $\|T-T'\|\leq\eps/3$.  Further, set $\phi_{n}(A) := \frac{Tr
    AP(VK_{n})} {Tr P(VK_{n})}$.  Then choose $n$ such that, for every
    $m>n$, $|\phi_{m}(T') - \phi_{n}(T')| \leq \eps/3$.  We get
    \begin{equation*}
	|\phi_{m}(T)-\phi_{n}(T)|\leq
	|\phi_{m}(T-T')| + |\phi_{m}(T')-\phi_{n}(T')| + 
	|\phi_{n}(T-T')|\leq\eps.
    \end{equation*}
    Hence, we have proved that $\lim \phi_{n}(T)$ exists.
    
    b) The functional $Tr_{\cg}$ is clearly linear, positive and
    takes value $1$ on the identity, hence it is a state on $\ca(X)$.
    \medskip
    
    \noindent {\bf Fourth step}: $Tr_{\cg}$ is a trace on $\ca(X)$. 
    
    Let $A$ be a geometric operator with propagation $r$. Then
    \begin{align*}
	AP(VK_{n}) &= P(B_{r}(VK_{n}))AP(VK_{n}), \\
	P(\O_{n,r})A &= P(\O_{n,r}) A P(VK_{n}).
    \end{align*}
    Indeed, the first equality is easily verified. As for the second, 
    we have
    $$
    \O_{n,r} \subset VK_{n} \setminus B_{r}(\cf(K_{n})) = \set{v\in
    VK_{n}: \r(v,VX\setminus VK_{n})\geq r+2},
    $$
    so that
    $$
    B_{r}(\O_{n,r}) \subset \set{v\in VK_{n}: \r(v,VX\setminus
    VK_{n})\geq 2} \subset VK_{n}.
    $$
    Since $A^{*}$ has propagation $r$, we get
    $$
    A^{*}P(\O_{n,r}) = P(B_{r}(\O_{n,r}))A^{*}P(\O_{n,r}) =
    P(VK_{n})A^{*}P(\O_{n,r}),
    $$
    which proves the claim.  Hence,
    \begin{align*}
	&AP(VK_{n}) = P(B_{r}(VK_{n}) \setminus \O_{n,r})AP(VK_{n}) +
	P(\O_{n,r})A \\
	&= 
	P(B_{r}(VK_{n}) \setminus \O_{n,r})AP(VK_{n})
	- P(VK_{n} \setminus \O_{n,r})A + P(VK_{n}) A.
    \end{align*}
    Therefore, if $B\in\ca (X)$,
    \begin{align*}
	\phi_{n}([B,A])&\leq\|A\|\,\|B\| \frac{|B_{r}(VK_{n})
	\setminus \O_{n,r}|+ |VK_{n} \setminus \O_{n,r}|}{|K_{n}|}\\
	&\leq 2 \|A\|\,\|B\|\eps_{n}(d+1)^{r},
    \end{align*}
    since $B_{r}(VK_{n}) \setminus \O_{n,r} \subset
    B_{r}(\cf_{\cg}(K_{n}))$.  Taking the limit as $n\to\infty$, we
    deduce that $Tr_{\cg}([B,A])=0$.  By continuity, it follows that
    the result holds for any $A,B\in\ca (X)$.
\end{proof}

\begin{rem}
	Let us note that the trace described in the previous theorem is not faithful in general. For example, for the Gasket graph in figure 1a the trace of $4I-D$ vanishes.
\end{rem}

\section[An analytic determinant for C$^*$-algebras]{An analytic
determinant for C$^*$-algebras with a trace state}

 In this section, we define a determinant for a suitable class of not
 necessarily normal operators in a C$^*$-algebra with a trace state, cf. \cite{Exel} for related results. 
 The results obtained are used in Section \ref{sec:detFormula} to
 prove a determinant formula for the zeta function.

 In a celebrated paper \cite{FuKa}, Fuglede and Kadison defined a
 positive-valued determinant for finite factors ($i.e.$ von Neumann
 algebras with trivial center and finite trace).  Such a determinant
 is defined on all invertible elements and enjoys the main properties
 of a determinant function, but it is positive-valued.  Indeed, for an
 invertible operator $A$ with polar decomposition $A=UH$, where $U$ is
 a unitary operator and $H:= \sqrt{A^{*}A}$ is a positive self-adjoint
 operator, the Fuglede--Kadison determinant is defined by
 $$
 Det(A)=\exp\, \circ\ \tau\circ\log H,
 $$
 where $\log H$ may be defined via the functional calculus.

 For the purposes of the present paper, we need a determinant which is
 an analytic function.  As we shall see, this can be achieved but
 corresponds to a restriction of the domain of the determinant
 function and implies the loss of some important properties.  Let
 $(\ca,\tau)$ be a C$^{*}$-algebra endowed with a trace state.  Then,
 a natural way to obtain an analytic function is to define, for
 $A\in\ca$, $\Det_\t(A)=\exp\, \circ\ \tau\circ\log A$, where
 $$
 \log(A) := \frac{1}{2\pi i} \int_\Gamma \log \lambda (\lambda-A)^{-1}
 d\lambda
 $$
 and $\Gamma$ is the boundary of a connected, simply connected region
 $\Omega$ containing the spectrum of $A$.  Clearly, once the branch of
 the logarithm is chosen, the integral above does not depend on
 $\Gamma$, provided $\G$ is given as above.

 Then a na\"{\i}ve way of defining $det$ is to allow all elements $A$
 for which there exists an $\Omega$ as above and a branch of the
 logarithm whose domain contains $\Omega$.  Indeed, the following
 holds.

\begin{Lemma}
    Let $A$, $\Omega$, $\Gamma$ be as above, and $\varphi$, $\psi$ two
    branches of the logarithm such that both associated domains
    contain $\Omega$.  Then
    $$
    \exp\, \circ\ \tau\circ\varphi(A) = \exp\, \circ\
    \tau\circ\psi(A).
    $$
\end{Lemma}

\begin{proof}
    The function $\varphi(\lambda)-\psi(\lambda)$ is continuous and
    everywhere defined on $\Gamma$.  Since it takes its values in
    $2\pi i \mathbb{Z}$, it should be constant on $\Gamma$. 
    Therefore,
    \begin{align*}
	\exp\, \circ\ \tau\circ\varphi (A) &= \exp\, \circ\
	\tau\left(\frac{1}{2\pi i} \int_\Gamma 2\pi i n_{0}
	(\lambda-A)^{-1} d\lambda \right) \exp\, \circ\
	\tau\circ\psi(A)\\
	&= \exp\, \circ\ \tau\circ\psi(A).
    \end{align*}
\end{proof}

 The problem with the previous definition is its dependence on the
 choice of $\Omega$.  Indeed, it is easy to see that when
 $A=\begin{pmatrix}1&0\\0&i\end{pmatrix}$, and if we choose $\Omega$
 containing $\{e^{i\vartheta},\vartheta\in[0,\pi/2]\}$ and any
 suitable branch of the logarithm, the determinant defined in terms of
 the normalized trace gives $det(A)=e^{i\pi/4}$.  On the other hand,
 if we choose $\Omega$ containing
 $\{e^{i\vartheta},\vartheta\in[\pi/2,2\pi]\}$ and a corresponding
 branch of the logarithm, we have $det(A)=e^{5i\pi/4}$.  Hence, we
 make the following choice.

\begin{Dfn}
    Let $(\ca,\tau)$ be a C$^{*}$-algebra endowed with a trace state,
    and consider the subset $\ca_{0}:=\{A\in\ca : 0\not\in
    \conv\sigma(A)\}$, where $\s(A)$ denotes the spectrum of $A$ and
    $\conv\s(A)$ its convex hull.  For any $A\in\ca_{0}$ we set
    $$
    \Det_\t(A)=\exp\, \circ\ \tau\circ\left(\frac{1}{2\pi i}
    \int_\Gamma \log \lambda (\lambda-A)^{-1} d\lambda\right),
    $$
    where $\Gamma$ is the boundary of a connected, simply connected
    region $\Omega$ containing $\conv\sigma(A)$, and $\log$ is a
    branch of the logarithm whose domain contains $\Omega$.
\end{Dfn}

Since two $\G$'s as above are homotopic in $\bc\setminus\conv\s(A)$,
we have

\begin{Cor}\label{cor:det.analytic}
    The determinant function defined above is well-defined and
    analytic on $\ca_{0}$.
\end{Cor}

 We collect several properties of our determinant in the following
 result.

\begin{Prop}
    Let $(\ca,\tau)$ be a C$^{*}$-algebra endowed with a trace state,
    and let $A\in\ca_{0}$.  Then
    
    \item[$(i)$] $\Det_\t(zA)=z\Det_\t(A)$, for any
    $z\in\mathbb{C}\setminus\{0\}$,
     
    \item[$(ii)$] if $A$ is normal, and $A=UH$ is its polar
    decomposition,
    $$
    \Det_\t(A)=\Det_\t(U)\Det_\t(H),
    $$
    
    \item[$(iii)$] if $A$ is positive, then we have
    $\Det_\t(A)=Det(A)$, where the latter is the Fuglede--Kadison
    determinant.
\end{Prop}

\begin{proof}
    $(i)$ If, for a given $\vartheta_0\in[0,2\pi)$, the half-line
    $\{\rho e^{i\vartheta_0}\in\mathbb{C} : \r>0\}$ does not intersect
    $\conv\sigma(A)$, then the half-line $\{\rho
    e^{i(\vartheta_0+t)}\in\mathbb{C} : \r>0\}$ does not intersect
    $\conv\sigma(zA)$, where $z=re^{it}$.  If $\log$ is the branch of
    the logarithm defined on the complement of the real negative
    half-line, then $\varphi(x)=i(\vartheta_{0}-\pi) +
    \log(e^{-i(\vartheta_{0}-\pi)}x)$ is suitable for defining
    $\Det_\t(A)$, while $\psi(x)=i(\vartheta_{0}+t-\pi) +
    \log(e^{-i(\vartheta_{0}+t-\pi)}x)$ is suitable for defining
    $\Det_\t(zA)$.  Moreover, if $\Gamma$ is the boundary of a
    connected, simply connected region $\Omega$ containing
    $\conv\sigma(A)$, then $z\Gamma$ is the boundary of a connected,
    simply connected region $z\Omega$ containing $\conv\sigma(zA)$. 
    Therefore,
    \begin{align*}
	\Det_\t(zA) &= \exp\, \circ\ \tau\left(\frac{1}{2\pi i}
	\int_{z\Gamma} \psi(\lambda) (\lambda-zA)^{-1}
	d\lambda\right)\\
	&= \exp\, \circ\ \tau\left(\frac{1}{2\pi i} \int_{\Gamma}
	(i(\vartheta_{0}+t-\pi) + \log(e^{-i(\vartheta_{0}+t-\pi)}
	re^{it}\mu)) (\mu-A)^{-1} d\mu\right)\\
	&= \exp\, \circ\ \tau\left((\log r + it)I+\frac{1}{2\pi i}
	\int_{\Gamma} \varphi(\mu) (\mu-A)^{-1} d\mu\right)\\
	&= z \Det_\t(A).
    \end{align*}
    $(ii)$ When $A=UH$ is normal, $U=\int_{[0,2\pi]} e^{i\vartheta}\
    du(\th)$, $H=\int_{[0,\infty)}r\ dh(r)$, then $ A =
    \int_{[0,\infty)\times[0,2\pi]} r e^{i\vartheta} \ d(h(r)\otimes
    u(\th))$.  The property $0\not\in\conv\sigma(A)$ is equivalent to
    the fact that the support of the measure $d(h(r)\otimes u(\th))$
    is compactly contained in some open half-plane 
    $$
    \{\rho e^{i\vartheta} : \rho>0, \vartheta \in (\vartheta_{0} -
    \pi/2, \vartheta_{0} +\pi/2)\},
    $$ 
    or, equivalently, that the support of the measure $dh(r)$ is
    compactly contained in $(0,\infty)$ and the support of the measure
    $d u(\th)$ is compactly contained in $(\vartheta_{0} - \pi/2,
    \vartheta_{0} +\pi/2)$.  Thus, $A\in\ca_{0}$ is equivalent to
    $U,H\in\ca_{0}$.  Then $$\log A = \int_{[0,\infty) \times
    (\vartheta_{0} - \pi/2, \vartheta_{0} +\pi/2)} (\log r +
    i\vartheta) \ d(h(r)\otimes u(\th)),$$ which implies that
    \begin{align*}
	\Det_\t(A) &= \exp\, \circ\ \tau\left(\int_{0}^{\infty} \log
	r\ dh(r) + \int_{\vartheta_{0} - \pi/2}^{\vartheta_{0} +\pi/2}
	i\vartheta \ du(\th)\right) \\
	&= \Det_\t(U)\cdot \Det_\t(H).
    \end{align*}
    $(iii)$  This follows by the  argument given in $(ii)$.
\end{proof}

\begin{rem} 
    We note that the above defined determinant function strongly
    violates the product property $\Det_\t(AB)=\Det_\t(A)\Det_\t(B)$. 
    Firstly, $A,B\in\ca_{0}$ does not imply $AB\in\ca_{0}$, as is seen
    e.g. by taking $A=B=\begin{pmatrix}1&0\\0&i\end{pmatrix}$. 
    Moreover, even if $A,B,AB\in\ca_{0}$ and $A$ and $B$ commute, the
    product property may be violated, as is shown by choosing
    $A=B=\begin{pmatrix}1&0\\0&e^{3i\pi/4}\end{pmatrix}$ and using
    the normalized trace on $2\times 2$ matrices.
\end{rem}
 
\section{Reduced closed paths}\label{sec:PrimeCycles}

 The Ihara zeta function is defined by means of equivalence classes of
 prime cycles.  Therefore, we need to introduce some terminology from
 graph theory, following \cite{StTe}, with some suitable
 modifications.  Moreover, we provide several technical results that
 will be used in the following sections.  We note that our approach to
 the determinant formula, here and in the following two sections,
 is inspired by the simplification brought to the subject by Stark and
 Terras in \cite{StTe}.
 
 Recall that a path (of length $m$) in $X$ from $v_0\in VX$ to $v_m\in
 VX$ is given by $(v_{0},\ldots,v_{m})$, where $v_{i}\in VX$,
 $v_{i+1}\sim v_i$, for $i=0,...,m-1$.  A path is closed if $v_m=v_0$. 
 Let us notice that an initial point $v_0$ is assigned also if the
 path is closed.
  
 \begin{Dfn}[Proper closed Paths] \label{def:redPath} 
    
    \itm{i} A path in $X$ has {\it backtracking} if $v_{i-1}=v_{i+1}$,
    for some $i\in\{1,\ldots,m-1\}$.  A path with no backtracking is
    also called {\it proper}.  Denote by $\Ci$ the set of proper
    closed paths, and by $\Ci_m$ the subset given by paths of length $m$.
       
    \itm{ii} A proper closed path $C=(v_{0},\ldots,v_{m}=v_{0})$ has a
    {\it tail} if there is $k\in\{1,\ldots,[\frac{m}{2}]-1\}$ such that
    $v_{j}=v_{m-j}$, for $j=1,\ldots,k$.  Denote by $\Ta$ the set of
    proper closed paths with tail, and by $\Nt$ the set of proper
    tail-less closed paths, also called {\it reduced closed paths}.
    The symbols $\Ta_m$ and $\Nt_m$ will denote the corresponding subsets 
    of paths of length $m$. 
    Observe that $\Ci=\Ta\cup\Nt$, $\Ta\cap\Nt=\emptyset$.
    
    \itm{iii} A reduced closed path is {\it primitive} if it is not
    obtained by going $n\geq 2$ times around some other closed path.
 \end{Dfn}

 Let us denote by $A$ the adjacency matrix of $X$, as defined in
 (\ref{eq:adjacency}).  For any $m\in\bn$, let us denote by
 $A_{m}(x,y)$ the number of proper paths in $X$, of length $m$, with
 initial vertex $x$ and terminal vertex $y$, for $x,y\in VX$.  Then
 $A_{1}=A$.  Let $A_{0}:= I$, and set $Q = D - I$, where $D$ is given
 by (\ref{eq:degree}).  Then, we obtain
 
 \begin{Lemma}\label{lem:Lemma1}
     \itm{i} $A_{2} = A^{2}-Q-I\in\ca(X)$,
     
     \itm{ii} for $m\geq 3$, $A_{m} = A_{m-1}A-A_{m-2}Q\in\ca(X)$,
     
     \itm{iii} if we let $d:= \sup_{v\in VX} \deg(v)$ and $\a:=
     \frac{d+\sqrt{d^2+4d}}{2}$, we have $\|A_{m}\| \leq \a^{m}$, for
     $m\geq0$.
 \end{Lemma}
 \begin{proof}
     $(i)$ If $x = y$ then $A_{2}(x,y)=0$ and $A^{2}(x,x) = deg(x) =
     (Q+I)(x,x)$.  If $x\neq y$, $A_{2}(x,y)=A(x,y)$ and
     $(Q+I)(x,y)=0$.
     
     $(ii)$ For $x,y\in VX$, the sum $\sum_{z\in VX}
     A_{m-1}(x,z)A(z,y)$ counts the proper paths of length $m$ from
     $x$ to $y$, which are $A_{m}(x,y)$, plus the paths of length $m$
     from $x$ to $y$ with backtracking at the end.  The paths in the
     latter family consist of a proper path of length $m-2$ from $x$
     to $y$ followed by a path of length $2$ from $y$ to $z$ and back;
     to avoid further backtracking, there are only $Q(y,y)$ choices
     for $z$, namely there are $A_{m-2}(x,y)Q(y,y)$ such paths.
     
     $(iii)$ We have $\|A_{1}\|=\|A\|\leq d$, $\|A_{2}\|\leq d^{2}+d$,
     and $\|A_{m}\|\leq d(\|A_{m-1}\|+\|A_{m-2}\|)$, from which the
     claim follows by induction.
 \end{proof}
 
 \begin{Lemma}\label{lem:countTail}
     For $m\in\bn$, let 
     $$
     t_{m}:= \lim_{n\to\infty} \frac{1}{|K_{n}|} \sum_{x\in K_{n}}
     \card{C\in\Ta_{m}: C \text{ starts at } x }.
     $$
     Then 
     
     \itm{i} the above limit exists and is finite,
    
     \itm{ii} $t_{1}=t_{2}=0$, and, for $m\geq 3$, $t_{m} = t_{m-2} +
     Tr_{\cg}((Q-I)A_{m-2})$,
     
     \itm{iii} $t_{m} = Tr_{\cg}\left(
     (Q-I)\sum_{j=1}^{[\frac{m-1}{2}]} A_{m-2j} \right)$, where 
     $Tr_{\cg}$ is defined in Theorem \ref{thm:trace}.
 \end{Lemma}
 \begin{proof}
     Denote by $(C,v)$ the proper closed path $C$ with the origin in
     $v\in VX$.
 
     $(i)$ For $n\in\bn$, let 
     $$
     \O_n:= V(K_n)\setminus B_1(\cf_\cg(K_n)), \qquad \O_{n}':=
     B_{1}(\cf_{\cg}(K_{n}))\cap V(K_{n}).
     $$ 
     Then, for all $p\in\bn$,
     $$
     V(K_{n+p}) = \left(\bigcup_{\g\in\cg(n,n+p)} \g\O_n\right) \cup
     \left(\bigcup_{\g\in\cg(n,n+p)} \g\O_{n}' \right).
     $$
     Let $t_m(x) := |\set{ (C,x)\in\Ta_m}|\leq d(d-1)^{m-2}$.  Then
     \begin{align*}
	& \left| \frac{1}{|K_{n+p}|} \sum_{x\in K_{n+p}} t_m(x) -
	\frac{1}{|K_{n}|} \sum_{x\in K_{n}} t_m(x) \right| \\
	& \leq \left| \frac{|\cg(n,n+p)|}{|K_{n+p}|} \sum_{x\in
	\O_{n}} t_m(x) - \frac{1}{|K_{n}|} \sum_{x\in K_{n}} t_m(x)
	\right| + \frac{|\cg(n,n+p)|}{|K_{n+p}|} \sum_{x\in
	\O_{n}'} t_m(x) \\
	& \leq \left| \frac{|\cg(n,n+p)|}{|K_{n+p}|} -
	\frac{1}{|K_{n}|} \right| \sum_{x\in K_{n}} t_m(x) + 2
	\frac{|\cg(n,n+p)|}{|K_{n+p}|} \sum_{x\in B_1(\cf_\cg(K_n))}
	t_m(x) \\
	& \leq \left| 1- \frac{|K_{n}| |\cg(n,n+p)|}{|K_{n+p}|}
	\right| d(d-1)^{m-2} + 2 \frac{|K_{n}| |\cg(n,n+p)|}{|K_{n+p}|}
	\frac{|B_1(\cf_\cg(K_n))|}{|K_n|} d(d-1)^{m-2} \\
	& \leq 6(d-1)^{m-2}d(d+1) \eps_n \to 0,\quad \text { as }
	n\to\infty,
    \end{align*}
    where, in the last inequality, we used equations (\ref{subsets}),
    (\ref{ineqE}) [with $r=1$], and the fact that $\eps_n =\displaystyle
    \frac{|\cf_\cg(K_n)|}{|K_n|}\to0$.
 
    $(ii)$ Let us define $\O:= \set{v\in VX: v\not\in K_n,
    \r(v,K_n)=1} \subset B_1(\cf_\cg(K_n))$.  We have
    \begin{align*}
	\frac{1}{|K_{n}|} \sum_{x\in K_{n}} \sum_{y\sim x}
	&|\set{C=(x,y,...)\in\Ta_{m} }| =\\
	&= \frac{1}{|K_{n}|} \sum_{y\in K_{n}} \sum_{x\sim y}
	|\set{C=(x,y,...)\in\Ta_{m} }| \\
	& \quad + \frac{1}{|K_{n}|} \sum_{y\in \O} \sum_{x\in K_n,
	x\sim y} |\set{C=(x,y,...)\in\Ta_{m} }| \\
	& \quad - \frac{1}{|K_{n}|} \sum_{y\in K_{n}} \sum_{x\in \O,
	x\sim y} |\set{C=(x,y,...)\in\Ta_{m} }| .
    \end{align*}
    Since 
    $$
    \frac{1}{|K_{n}|} \sum_{y\in\O} \sum_{x\in K_n, x\sim y}
    |\set{C=(x,y,...)\in\Ta_{m} }| \leq \frac{1}{|K_n|}
    |\cf_\cg(K_n)|d(d+1)(d-1)^{m-2} \to0
    $$
    and
    \begin{align*}
	 \frac{1}{|K_{n}|} \sum_{y\in K_{n}} \sum_{x\in \O, x\sim y}
	&|\set{C=(x,y,...)\in\Ta_{m} }| = \\
	& = \frac{1}{|K_{n}|} \sum_{y\in \cf_\cg(K_{n})} \sum_{x\in
	\O, x\sim y} |\set{C=(x,y,...)\in\Ta_{m} }| \\
	& \leq \frac{1}{|K_n|} |\cf_    \cg(K_n)|d(d-1)^{m-2} \to0,
   \end{align*}
   we obtain
   \begin{align*}
       t_{m} &= \lim_{n\to\infty} \frac{1}{|K_{n}|} \sum_{x\in K_{n}}
       |\set{ (C,x)\in\Ta_m }| \\
       &= \lim_{n\to\infty} \frac{1}{|K_{n}|} \sum_{x\in K_{n}}
       \sum_{y\sim x} |\set{C=(x,y,...)\in\Ta_{m} }| \\
       &= \lim_{n\to\infty} \frac{1}{|K_{n}|} \sum_{y\in K_{n}}
       \sum_{x\sim y} |\set{C=(x,y,...)\in\Ta_{m} }| .
   \end{align*}

   For a given $y$, we now want to count the paths in $\Ta_m$ whose
   second vertex is $y$.  Any such path can be identified by the
   choice of a first vertex $x\sim y$ and a proper closed path $D$ of
   length $(m-2)$\footnote{Here and thereafter, $D$ is a path of the
   graph and not the degree matrix.}.  Of course we may first choose
   the closed path $D$ and then the first vertex $x$.  There are two
   kinds of proper closed paths $D$ starting at $y$, namely, those
   with tails and those without.  If $D$ has no tail, then there are
   $Q(y,y)-1$ possibilities for $x$ to be adjacent to $y$ in such a
   way that the resulting path $C$ has no backtracking.  If $D$ has a
   tail, then there are $Q(y,y)$ possibilities for $x$ to be adjacent
   to $y$ in such a way that the resulting path $C$ has no
   backtracking.  Therefore, we have
   \begin{align*}
       \sum_{x\sim y} &|\set{C=(x,y,...)\in\Ta_{m} }| = \\
       &= (Q(y,y)-1)\cdot |\set{(D,y)\in\Nt_{m-2} }| + Q(y,y)\cdot
       |\set{(D,y)\in\Ta_{m-2} }| \\
       & = (Q(y,y)-1)\cdot |\set{ (D,y)\in\Ci_{m-2} }| +
       |\set{(D,y)\in\Ta_{m-2}} | \, ,
   \end{align*}
   so that
   \begin{align*}
       t_{m} &= \lim_{n\to\infty} \frac{1}{|K_{n}|} \sum_{y\in K_{n}}
       (Q(y,y)-1)\cdot |\set{ (D,y)\in\Ci_{m-2}}| \\
       & \qquad + \lim_{n\to\infty} \frac{1}{|K_{n}|} \sum_{y\in
       K_{n}} |\set{(D,y)\in\Ta_{m-2}}| \\
       & = \lim_{n\to\infty} \frac{1}{|K_{n}|} \sum_{y\in K_{n}}
       (Q(y,y)-1)A_{m-2}(y,y) +t_{m-2}\\
       & = Tr_{\cg}((Q-I)A_{m-2})+t_{m-2}.
   \end{align*}
   $(iii)$ This follows from $(ii)$.
\end{proof}
 
 \begin{Lemma}\label{lem:estim.for.N}
     Let 
     $$
     N_{m} := \lim_{n\to\infty} \frac{1}{|K_{n}|} | \set{ C\in\Nt_{m}:
     C\subset K_{n} }|,
     $$
     which exists and is finite.  Then, for all $m\in\bn$,
     
     \itm{i} $N_{m} = Tr_{\cg}(A_{m}) - t_{m}$,
     
     \itm{ii} $N_{m} \leq d(d-1)^{m-1}$.
 \end{Lemma}
 \begin{proof}
     $(i)$ Observe that 
     \begin{equation}\label{eq:Nm}
	 N_{m} = \lim_{n\to\infty} \frac{1}{|K_{n}|} |\set{
	 (C,v)\in\Nt_{m}: v\in K_{n} }|.
     \end{equation}
     Indeed, 
     \begin{align*}
	 0&\leq \frac{1}{|K_{n}|} \Big| |\set{ (C,v)\in\Nt_{m}: v\in
	 K_{n} }| - | \set{ C\in\Nt_{m}: C\subset K_{n} }|  \Big| \\
	 & = \frac{1}{|K_{n}|}\, \card{ (C,v)\in\Nt_{m}: v\in
	 K_{n}, C\not\subset K_{n} }\\
	 & \leq \frac{1}{|K_{n}|}\, \card{ (C,v)\in\Ci_{m}: v\in
	 B_m(\cf_\cg(K_n)) }\\
	 & = \frac{1}{|K_{n}|}\, \sum_{v\in B_m(\cf_\cg(K_n))}
	 A_{m}(v,v) = \frac{1}{|K_{n}|}\, Tr(P( B_m(\cf_\cg(K_n)))
	 A_{m}) \\
	 & \leq \|A_{m}\|\, \frac{|B_m(\cf_\cg(K_n))|}{|K_{n}|} 
	 \leq \|A_{m}\|\, (d+1)^m \frac{| \cf_\cg(K_n)|}{|K_{n}|}\to
	 0,\, \text{as } n\to\infty.
     \end{align*}
     Moreover, the existence of $\lim_{n\to\infty} \frac{1}{|K_{n}|}
     |\set{ (C,v)\in\Nt_{m}: v\in K_{n} }|$ can be proved as in Lemma
     \ref{lem:countTail} $(i)$.  Therefore,
     \begin{align*}
	 N_{m} & = \lim_{n\to\infty} \frac{1}{|K_{n}|} \card{
	 (C,v)\in\Nt_{m}: v\in K_{n} } \\
	 & = \lim_{n\to\infty} \frac{1}{|K_{n}|}\, \biggl( \sum_{v\in
	 K_{n}} A_{m}(v,v) - \sum_{v\in K_{n}} \card{ (C,v)\in\Ta_{m}:
	 v\in K_{n} } \biggr)\\
	 & = Tr_{\cg}(A_{m}) - t_{m}.	 
     \end{align*}
     $(ii)$ This follows from (\ref{eq:Nm}). 
 \end{proof}

\section{The Zeta function}\label{sec:Zeta}

 In this section, we define the Ihara zeta function for a self-similar
 graph and prove that it is a holomorphic function in a suitable disc. 
 We first need to introduce some equivalence relations between closed
 paths.
 
  \begin{Dfn}[Cycles]
    We say that two closed paths $C=(v_{0},\ldots,v_{m}=v_{0})$ and
    $D=(w_{0},\ldots,w_{m}=w_{0})$ are {\it equivalent}, and write
    $C\sim_{o} D$, if there is an integer $k$ such that
    $w_{j}=v_{j+k}$, for all $j$, where the addition is taken modulo
    $m$, that is, the origin of $D$ is shifted $k$ steps with respect
    to the origin of $C$.  The equivalence class of $C$ is denoted
    $[C]_o$.  An equivalence class is also called a {\it cycle}.
    Therefore, a closed path is just a cycle with a specified origin.

    Denote by $\Re$ the set of reduced cycles, and by $\Pr\subset\Re$
    the subset of primitive reduced cycles, also called {\it prime}
    cycles.
 \end{Dfn}

 \begin{Dfn}[Equivalence relation between reduced cycles]
    Given $C$, $D\in\Re$, we say that $C$ and $D$ are $\cg$-{\it
    equivalent}, and write $C \sim_{\cg} D$, if there is a local
    isomorphism $\g\in \cg$ such that $D=\g(C)$.  We denote by
    $[\Re]_{\cg}$ the set of $\cg$-equivalence classes of reduced
    cycles, and analogously for the subset $\Pr$.
 \end{Dfn}
 
  We need to introduce several quantities associated to a reduced cycle.
 
  \begin{Dfn}[Average multiplicity of a reduced cycle]
     Let $C\in\Re$, and call 
     
     \itm{i} {\it size} of $C$, denoted $s(C)\in\bn$, the least $m\in\bn$
     such that $C\subset \g(K_{m})$, for some local isomorphism $\g\in
     \cg(m)$,
     
     \itm{ii} {\it effective length} of $C$, denoted $\ell(C)\in\bn$,
     the length of the prime cycle $D$ underlying $C$, $i.e.$ such
     that $C=D^p$, for some $p\in\bn$,
     
     \itm{iii} {\it average multiplicity} of $C$, denoted $\mu(C)$,
     the number in $[0,\infty)$ given by
     $$
     \lim_{n\to\infty} \frac{|\cg(s(C),n)|}{|K_{n}|}.
     $$
 \end{Dfn}

That the limit actually exists is the content of the following
 
 \begin{Prop}
     \itm{i} Let $C\in\Re$, then the following limit exists and is
     finite:
     $$
     \lim_{n}\frac{ |\cg(s(C),n)| }{|K_n|},
     $$

     \itm{ii} $s(C)$, $\ell(C)$, and $\m(C)$ only depend on
     $[C]_{\cg}\in [\Re]_{\cg}$; moreover, if $C=D^{k}$ for some
     $D\in\Pr$, $k\in\bn$, then $s(C)=s(D)$, $\ell(C)=\ell(D)$,
     $\m(C)=\m(D)$,
         
     \itm{iii} for $m\in\bn$, $\displaystyle{N_m =
     \sum_{[C]_\cg\in[\Re_m]_\cg} \m(C)\ell(C)}$,\\
     where, as above, the subscript $m$ corresponds to reduced cycles
     of length $m$.
 \end{Prop}
 \begin{proof}
     $(i)$ Let us observe that $|\cg(s(C),n+1)| = |\cg(s(C),n)|
     |\cg(n,n+1)|$, for any integer $n\geq s(C)$.  Therefore, using
     (\ref{ineqE}) and (\ref{ineqC}), we obtain
     \begin{align}\label{estimate}
	 \left| \frac{|\cg(s(C),n)|} {|K_n|} - \frac{|\cg(s(C),n+p)|}
	 {|K_{n+p} |} \right| & = \frac{|\cg(s(C),n)|} {|K_n|} \left|
	 1- \frac{ |K_n| |\cg(n,n+p)| }{|K_{n+p}|} \right| \notag\\
     \leq \frac{1}{| \O_{s(C),1} |} 2\eps_n (d+1).
     \end{align}
     The regular exhaustion property implies the claim.  Let us
     finally observe that the limit is monotone.  Indeed, for $n\geq
     s(C)$,
     \begin{align*}
	 \frac{|\cg(s(C),n+1)|} {|K_{n+1}|} = \frac{|\cg(s(C),n)|}
	 {|K_{n} |} \frac{ |K_n| |\cg(n,n+1)| }{|K_{n+1}|} \geq
	 \frac{|\cg(s(C),n)|} {|K_{n} |} .
     \end{align*}
     
     $(ii)$ This follows from the definition.
     
     $(iii)$ We have successively:
     \begin{align*} 
	 N_m & = \lim_{n\to\infty} \frac{1}{|K_{n}|} | \set{
	 C\in\Nt_{m}: C\subset K_{n} }| \\
	 & = \lim_{n\to\infty} \sum_{[C]_\cg\in[\Re_m]_\cg}
	 \frac{1}{|K_{n}|} | \set{ D\in\Nt_{m}: D\sim_\cg C, D\subset
	 K_{n} }| \\
	 & = \lim_{n\to\infty} \sum_{[C]_\cg\in[\Re_m]_\cg}
	 \frac{1}{|K_{n}|} \, \ell(C)\, | \cg(s(C),n)|  \\
	 & = \sum_{[C]_\cg\in[\Re_m]_\cg} \m(C)\ell(C),
     \end{align*} 
     where, in the last equality, we used monotone convergence.  
  \end{proof}
 
 \begin{exmp}
     $(i)$ For the Gasket graph of figure \ref{fig:GasketVicsek}a we
     get, for a prime cycle $C$ of size $s(C)=p$, that $\displaystyle{
     \mu(C) = \lim_{n\to\infty} \frac{3^{n-p}} {\frac32 (3^{n}+1)} =
     \frac{2}{3^{p+1}} }$.  \\
     $(ii)$ For the Vicsek graph of figure \ref{fig:GasketVicsek}b we
     get, for a prime cycle $C$ of size $s(C)=p$, that $\displaystyle{
     \mu(C) = \lim_{n\to\infty} \frac{5^{n-p}} {3\cdot 5^{n}+1} =
     \frac{1}{3\cdot 5^{p}} }$.  \\
     $(iii)$ For the Lindstrom graph of figure
     \ref{fig:LindstromCarpet}a we get, for a prime cycle $C$ of size
     $s(C)=p$, that $\displaystyle{ \mu(C) = \lim_{n\to\infty}
     \frac{7^{n-p}} {5\cdot 7^{n}+1} = \frac{1}{5\cdot 7^{p}} }$.
 \end{exmp}

 We can now introduce the counterpart of the Ihara zeta function for a
 self-similar graph.

 \begin{Dfn}[Zeta function]\label{Dfn:Zeta}
     Let $Z(u)=Z_{X,\cg}(u)$ be given by
     $$
     Z_{X,\cg}(u) := \prod_{[C]_{\cg}\in [\Pr]_{\cg}}
     (1-u^{|C|})^{-\mu(C)}, 
     $$
     for $u\in\bc$ sufficiently small so that the infinite product
     converges.
 \end{Dfn}

 \begin{Thm}\label{lem:power.series}
     Let $X$ be a self-similar graph, with $d:=\sup_{v\in VX}
     \deg(v)$.  Then, \itm{i} $Z(u):=\prod_{[C]_{\cg}\in [\Pr]_{\cg}}
     (1-u^{|C|})^{-\mu(C)}$ defines a holomorphic function in the open
     disc $\{u\in\bc: |u|<\frac{1}{d-1}\}$,
     
     \itm{ii} $u\frac{Z'(u)}{Z(u)} = \sum_{m=1}^{\infty} N_{m}u^{m}$,
     for $|u|<\frac{1}{d-1}$,
          
     \itm{iii} $Z(u) = \exp\left( \sum_{m=1}^{\infty} \frac{N_{m}}{m}
     u^{m} \right)$, for $|u|<\frac{1}{d-1}$.
 \end{Thm}
 \begin{proof}
     Let us observe that, for $|u|<\frac{1}{d-1}$,
     \begin{align*}
	 \sum_{m=1}^{\infty} N_{m} u^{m} & = \sum_{m=1}^{\infty}
	 \sum_{[C]_{\cg}\in [\Re_{m}]_{\cg}} \ell(C) \mu(C)\, u^{m} \\
	 & = \sum_{[C]_{\cg}\in [\Re]_{\cg}} \ell(C) \mu(C)\, u^{|C|} \\
	 & = \sum_{m=1}^{\infty} \sum_{[C]_{\cg}\in [\Pr]_{\cg}} 
	 |C| \mu(C)\,	 u^{|C^{m}|} \\
	 &= \sum_{[C]_{\cg}\in [\Pr]_{\cg}} \mu(C)\, \sum_{m=1}^{\infty}
	 |C| u^{|C|m} \\
	 &= \sum_{[C]_{\cg}\in [\Pr]_{\cg}} \mu(C)\, u\frac{d}{du}
	 \sum_{m=1}^{\infty} \frac{u^{|C|m}}{m} \\
	 &= -\sum_{[C]_{\cg}\in [\Pr]_{\cg}} \mu(C)\, u\frac{d}{du}
	 \log(1-u^{|C|}) \\
	 & = u\frac{d}{du} \log Z(u),
     \end{align*}
     where we have used the Fubini--Tonelli theorem in the fourth
     equality, while, in the last equality, we have used uniform
     convergence on compact subsets of $\set{u\in\bc:
     |u|<\frac{1}{d-1}}$.  The rest of the proof is clear.
 \end{proof}

\section{The determinant formula}\label{sec:detFormula}

 In this section, we establish the main result in the theory of Ihara
 zeta functions, which states that $Z$ is the reciprocal of a
 holomorphic function which, up to a multiplicative factor, is the
 determinant of a deformed Laplacian on the graph.  We first need to
 state several technical results.  Let us recall that $d:=\sup_{v\in
 VX} \deg(v)$ and $\a:= \frac{d+\sqrt{d^{2}+4d}}{2}$.
 
 \begin{Lemma}\label{lem:eq.for.A}
     \itm{i} $\left(\sum_{m\geq 0} A_{m}u^{m}\right)(I-Au+Qu^{2}) = 
     (1-u^{2})I$, $|u|<\frac{1}{\a}$, 
          
     \itm{ii} $\left(\sum_{m\geq 0} \left( 
     \sum_{k=0}^{[m/2]}A_{m-2k} \right) u^{m}\right)(I-Au+Qu^{2}) = 
     I$, $|u|<\frac{1}{\a}$.
 \end{Lemma}
 \begin{proof}
     $(i)$ From Lemma \ref{lem:Lemma1} we obtain
     \begin{align*}
	 \biggl(\sum_{m\geq 0} A_{m}u^{m}\biggr)&(I-Au+Qu^{2}) =
	 \sum_{m\geq 0} A_{m}u^{m} - \sum_{m\geq 0}\left(
	 A_{m}Au^{m+1} -A_{m}Qu^{m+2}\right) \\
	 &= \sum_{m\geq 0} A_{m}u^{m} -A_{0}Au -A_{1}Au^{2}
	 +A_{0}Qu^{2} \\
	 & \qquad - \sum_{m\geq 3}\left( A_{m-1}A
	 -A_{m-2}Q\right)u^{m} \\
	 &= \sum_{m\geq 0} A_{m}u^{m} -Au -A^{2}u^{2}
	 +Qu^{2} - \sum_{m\geq 3} A_{m}u^{m} \\
	 &= I +Au +A_{2}u^{2} -Au -A^{2}u^{2} +Qu^{2} \\
	 & = (1-u^{2})I.
     \end{align*}
     $(ii)$ We have successively:
     \begin{align*}
	 I &= (1-u^{2})^{-1} \biggl(\sum_{m\geq 0} 
	 A_{m}u^{m}\biggr)(I-Au+Qu^{2}) \\
	 &= \biggl(\sum_{m\geq 0} A_{m}u^{m}\biggr) \biggl(
	 \sum_{j=0}^{\infty}u^{2j}\biggr) (I-Au+Qu^{2}) \\
	 &= \biggl(\sum_{k\geq 0}\sum_{j=0}^{\infty} 
	 A_{k}u^{k+2j}\biggr)(I-Au+Qu^{2}) \\
	 &= \biggl(\sum_{m\geq 0}\biggl( \sum_{j=0}^{[m/2]} 
	 A_{m-2j}\biggr) u^{m}\biggr)(I-Au+Qu^{2}).
     \end{align*}     
 \end{proof}
 
 \begin{Lemma}\label{lem:eq.for.B}
     Let $B_{m} := A_{m} - (Q-I) \sum_{k=1}^{[m/2]}A_{m-2k}
     \in\ca(X)$, for $m\geq 0$.  Then 
     
     \itm{i} $B_{0}=I$, $B_{1}=A$,
     
     \itm{ii} $B_{m} = QA_{m} - (Q-I) \sum_{k=0}^{[m/2]}A_{m-2k}$,
     
     \itm{iii} $$Tr_{\cg} B_{m} = 
     \begin{cases}
	 N_{m} - Tr_{\cg}(Q-I) & m \text{ even} \\
	 N_{m}  & m \text{ odd,}
     \end{cases}$$
     \itm{iv}
     $$
     \sum_{m\geq 1} B_{m}u^{m} = \left(
     Au-2Qu^{2}\right)\left(I-Au+Qu^{2}\right)^{-1}, \
     |u|<\frac{1}{\a}.
     $$
 \end{Lemma}
 \begin{proof}
     $(i)$ and $(ii)$ follow from straightforward computations
     involving bounded operators.
     
     $(iii)$ It follows from Lemma \ref{lem:countTail} that, if $m$ is
     odd, 
     $$
     Tr_{\cg} B_{m} = Tr_{\cg}(A_{m}) - t_{m} = N_{m},
     $$
     whereas, if $m$ is even, 
     $$
     Tr_{\cg} B_{m} = Tr_{\cg}(A_{m}) - t_{m} - Tr_{\cg}((Q-I)A_{0}) =
     N_{m} - Tr_{\cg}(Q-I).
     $$
     $(iv)$ 
     \begin{align*}
	 \biggl( \sum_{m\geq 0} B_{m}u^{m} \biggr)& (I-Au+Qu^{2}) 
	 \\
	 & = \biggl( Q\sum_{m\geq 0} A_{m}u^{m} - (Q-I)\sum_{m\geq
	 0}\sum_{j=0}^{[m/2]} A_{m-2j}u^{m}\biggr) (I-Au+Qu^{2}) \\
	 & = Q(1-u^{2})I - (Q-I)\biggl( \sum_{m\geq 0}\sum_{j=0}^{[m/2]}
	 A_{m-2j}u^{m}\biggr) (I-Au+Qu^{2}) \\
	 & = (1-u^{2})Q - (Q-I) = I-u^{2}Q,
     \end{align*}
     where we used Lemma \ref{lem:eq.for.A} $(i)$, in the second
     equality, and Lemma \ref{lem:eq.for.A} $(ii)$, in the third
     equality.  Since $B_{0}=I$, we deduce that
     \begin{align*}
	 \biggl( \sum_{m\geq 1} B_{m}u^{m} \biggr) (I-Au+Qu^{2}) &=
	 I-u^{2}Q - B_{0}(I-Au+Qu^{2})\\
	 & = Au-2Qu^{2}.
     \end{align*}
 \end{proof}
 
 \begin{Lemma}\label{lem:Lemma3}
     Let $f:u\in B_{\eps} := \{u\in\bc: |u|<\eps\} \mapsto
     f(u)\in\ca(X)$, be a $C^{1}$-~function such that $f(0)=0$ and
     $\|f(u)\|<1$, for all $u\in B_{\eps}$.  Then
     $$
     Tr_{\cg}\left(
     -\frac{d}{du} \log(I-f(u)) \right) = Tr_{\cg}\left(
     f'(u)(I-f(u))^{-1}\right).
     $$
 \end{Lemma}
 \begin{proof}
     To begin with, $-\log(I-f(u)) = \sum_{n\geq 1} \frac{1}{n}
     f(u)^{n}$ converges in operator norm, uniformly on compact
     subsets of $B_{\eps}$.  Moreover, 
     $$
     \frac{d}{du} f(u)^{n} =
     \sum_{j=0}^{n-1} f(u)^{j}f'(u) f(u)^{n-j-1}.
     $$  
     Therefore,
     $-\frac{d}{du} \log(I-f(u)) = \sum_{n\geq 1} \frac{1}{n}
     \sum_{j=0}^{n-1} f(u)^{j}f'(u) f(u)^{n-j-1}$, so that
     \begin{align*}
	 Tr_{\cg}\biggl( -\frac{d}{du} \log(I-f(u)) \biggr) & =
	 \sum_{n\geq 1} \frac{1}{n} \sum_{j=0}^{n-1} Tr_{\cg}\left(
	 f(u)^{j}f'(u) f(u)^{n-j-1} \right) \\
	 & = \sum_{n\geq 1} Tr_{\cg}( f(u)^{n-1}f'(u) ) \\
	 & = Tr_{\cg}\biggl( \sum_{n\geq 0}  f(u)^{n}f'(u) \biggr) \\
	 & = Tr_{\cg}( f'(u)(I-f(u))^{-1} ),
     \end{align*}
     where we have used the fact that $Tr_{\cg}$ is norm continuous.
 \end{proof}
 
 \begin{Cor} 
     $$
     Tr_{\cg}\left( \sum_{m\geq 1} B_{m}u^{m} \right) = Tr_{\cg}\left(
     -u\frac{d}{du} \log(I-Au+Qu^{2}) \right),\ |u|<\frac{1}{\a}.
     $$
 \end{Cor}
 \begin{proof}
     It follows from Lemma \ref{lem:eq.for.B} $(iv)$ that
     \begin{align*}
	 Tr_{\cg}\biggl( \sum_{m\geq 1} B_{m}u^{m} \biggr) &=
	 Tr_{\cg}( (Au-2Qu^{2}) (I-Au+Qu^{2})^{-1} ) \\
	 &= Tr_{\cg} \Bigl( -u\frac{d}{du} \log(I-Au+Qu^{2}) \Bigr),
     \end{align*} 
     where the last equality follows from the previous lemma applied with
     $f(u) := Au-Qu^{2}$.
 \end{proof}
 
 We now introduce the average Euler--Poincar\'e characteristic of a
 self-similar graph.
 
 \begin{Lemma} \label{lem:EuPoi}
     The following limit exists and is finite: 
     $$
     \chi_{av}(X) := \lim_{n\to\infty} \frac{\chi(K_n)}{|K_n|} =
     -\frac12 Tr_\cg(Q-I),
     $$
     where $\chi(K_n)=|VK_n| - |EK_n|$ is the Euler--Poincar\'e
     characteristic of the subgraph $K_n$.  The number $\chi_{av}(X)$
     is called the average Euler--Poincar\'e characteristic of the
     self-similar graph $X$.
 \end{Lemma}
 \begin{proof}
     Let, for $v,\ w\in VK_{n}$, 
     $$
     Q_n(v,w) :=
     \begin{cases}
	 \deg(v)-1 & v=w\\
	 0& v\neq w,
     \end{cases}
     $$
     and let $\d_n:= (Q-Q_n)P(\cf_\cg K_n)$.  Hence,
     $QP(K_n\setminus \cf_\cg K_n) = Q_nP(K_n\setminus \cf_\cg K_n)$,
     and $QP(K_n) = Q_n +\d_n$.  Since 
     $$
     Tr(Q_n) = \sum_{v\in K_n} \deg(v) - |VK_n| = 2|EK_n| -|VK_n|,
     $$
     and 
     $$
     |Tr(\d_n)| \leq \|Q-Q_n\| Tr(P( \cf_\cg K_n)) \leq (d-1) |\cf_\cg
     K_n|,
     $$
     so that
     $ \lim_{n\to\infty} \frac{Tr(\d_n)}{|VK_n|} =0$, we obtain
     \begin{align*}
	 \lim_{n\to\infty} \frac{\chi(K_n)}{|VK_n|} & =
	 \lim_{n\to\infty} \frac{|VK_n|- |EK_n|}{|VK_n|} 
	 = \frac12  - \frac12 \lim_{n\to\infty} \frac{Tr(Q_n)}{|VK_n|} \\
	 & = \frac12 - \frac12 \lim_{n\to\infty}
	 \frac{Tr(Q_n+\d_n)}{|VK_n|} = \frac12 - \frac12
	 \lim_{n\to\infty} \frac{Tr(QP(K_n))}{|VK_n|} \\
	 & = - \frac12 Tr_\cg(Q-I).
     \end{align*}
 \end{proof}
 
 \begin{exmp}
     We compute the average Euler--Poincar\'e characteristic of some
     self-similar graphs.
 
     \itm{i} For the Gasket graph of figure \ref{fig:GasketVicsek}a,
     we get $|VK_{n}| = \frac12 3^{n}+\frac32$ and $|EK_{n}| =3^{n}$,
     so that $\chi_{av}(X)=-1$, see also Proposition
     \ref{prop:essentiallyRegular} $(i)$ and Example
     \ref{exmp:chi(gasket)}.
     
     \itm{ii} For the Vicsek graph of figure \ref{fig:GasketVicsek}b,
     we get $|VK_{n}| = 3\cdot 5^{n}+1$ and $|EK_{n}| =4\cdot 5^{n}$,
     so that $\chi_{av}(X)=-\frac13$.
     
     \itm{iii} For the Lindstrom graph of figure
     \ref{fig:LindstromCarpet}a, we get $|VK_{n}| = 4\cdot 7^{n}+2$
     and $|EK_{n}| =6\cdot 7^{n}$, so that $\chi_{av}(X)=-\frac12$.
     
     \itm{iv} For the Carpet graph of figure
     \ref{fig:LindstromCarpet}b, we get $|VK_{n}| = \frac{44}{35}
     8^{n}(1+o(1))$ and $|EK_{n}| =\frac{12}{5}8^{n}(1+o(1))$, so that
     $\chi_{av}(X)=-\frac{10}{11}$.  Here, as usual, $o(1)$ denotes a
     sequence tending to zero as $n\to\infty$.
 \end{exmp}
 
 \begin{rem}
     We note that the average Euler--Poincar\'e characteristic
     introduced in the previous Lemma \ref{lem:EuPoi} coincides with
     the L$^2$-Euler--Poincar\'e characteristic defined as the
     alternating sum of the L$^2$-Betti numbers, as shown in
     \cite{CGIs01}, though with a different normalization factor.
 \end{rem}

Recall that $d:=\sup_{v\in VX} \deg(v)$ and $\a:=
\frac{d+\sqrt{d^{2}+4d}}{2}$.
  
 \begin{Thm}[Determinant formula]
     Let $X$ be a self-similar graph and $Z_{X,\cg}$ its zeta
     function.  Then,
     $$
     \frac{1}{Z_{X,\cg}(u)} = (1-u^{2})^{-\chi_{av}(X)}
     \Det_{\cg}(I-Au+Qu^{2}), \ \text{ for }|u|<\frac{1}{\a}.
     $$ 
 \end{Thm}
 \begin{proof}
     We have
     \begin{align*}
	 Tr_{\cg}\biggl( \sum_{m\geq 1} B_{m}u^{m} \biggr) &=
	 \sum_{m\geq 1} Tr_{\cg}( B_{m} ) u^{m}= \sum_{m\geq 1}
	 N_{m}u^{m} - \sum_{k\geq 1} Tr_{\cg}(Q-I) u^{2k} \\
	 &= \sum_{m\geq 1} N_{m}u^{m} - Tr_{\cg}(Q-I)
	 \frac{u^{2}}{1-u^{2}},
     \end{align*}
     where the second equality follows from Lemma \ref{lem:eq.for.B}
     $(iii)$.  Therefore,
     \begin{align*}
	 u\frac{d}{du} &\log Z_{X,\cg}(u)  = \sum_{m\geq 1} N_{m}u^{m} \\
	 &= Tr_{\cg}\left( -u\frac{d}{du} \log(I-Au+Qu^{2}) \right) -
	 \frac{u}{2}\frac{d}{du} \log(1-u^{2}) Tr_{\cg}(Q-I)
     \end{align*}
     so that, dividing by $u$ and integrating from $u=0$ to $u$, we
     obtain
     $$
     \log Z_{X,\cg}(u) = - Tr_{\cg}\left( \log(I-Au+Qu^{2}) \right)
     -\frac12 Tr_{\cg}(Q-I) \log(1-u^{2})
     $$
     which implies that
     $$
     \frac{1}{Z_{X,\cg}(u)} = (1-u^{2})^{\frac12 Tr_{\cg}(Q-I)}
     \cdot\exp Tr_{\cg} \log(I-Au+Qu^{2}).
     $$
 \end{proof}
 
 \begin{rem}
     Observe that the domain of validity of the determinant formula
     above is smaller than the domain of holomorphicity of
     $Z_{X,\cg}$, which is the open disc $\{u\in\bc:
     |u|<\frac{1}{d-1}\}$, as was proved in Theorem
     \ref{lem:power.series}.
 \end{rem}
  
\section{Essentially regular graphs}\label{sec:essRegGraph}

 In this section, we obtain several functional equations for the Ihara
 zeta functions of essentially $(q+1)$-regular graphs, $i.e.$
 self-similar graphs $X$ such that $deg(v)=q+1$ for all but a finite
 number of vertices $v\in VX$.  The various functional equations
 correspond to different ways of completing the zeta functions.
  
\begin{Lemma} \label{prop:holomorphy}
    Let $X$ be essentially $(q+1)$-regular.  Then 
    $$
    u\in\O\mapsto \Det_{\cg}((1+qu^{2})I-Au)\in\bc
    $$ 
    is a holomorphic function at
    least in the open set
    $$
    \O:=\br^2 \setminus \left(\set{(x,y)\in\br^2: x^2+y^2=\frac{1}{q}}
    \cup \set{(x,0)\in\br^2: \frac{1}{q}\leq |x|\leq 1}\right).
    $$ 
    See figure \ref{fig:Omega}.
     \begin{figure}[ht]
	 \centering
	 \psfig{file=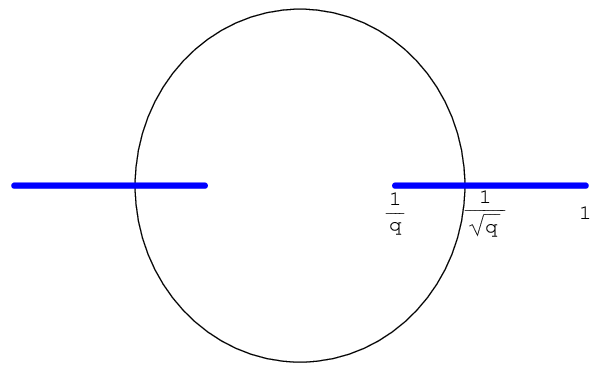,height=1.5in}
	 \caption{The open set $\O$ in Lemma \ref{prop:holomorphy}}
	 \label{fig:Omega}
     \end{figure}
\end{Lemma}
\begin{proof}
    Let $\D(u) := (1+qu^2)I-Au$, and observe that 
    $$
    \s(\D(u)) =
    \set{1+qu^2-u\l: \l\in\s(A)} \subset \set{1+qu^2-u\l:
    \l\in[-d,d]},
    $$
    so that $0\not\in\conv\s(\D(u))$ at least for $u\in\bc$ such that
    $1+qu^2-u\l\neq0$ for $\l\in[-d,d]$, that is for $u=0$ or
    $\frac{1+qu^2}{u}\not\in[-d,d]$, and hence at least for $u\in\O$. 
    The rest of the proof follows from Corollary
    \ref{cor:det.analytic}.
\end{proof}

 Let us denote by $P$ the {\it transition probability operator} of the
 simple random walk on $X$, that is
 $$
 P(x,y):= 
 \begin{cases}
    \frac{1}{\deg(x)} & y \text{ adjacent to } x\\
    0 & \text{otherwise.}
 \end{cases}
 $$

\begin{Prop} \label{prop:essentiallyRegular}
    Let $X$ be essentially $(q+1)$-regular, $i.e.$ $deg(v)=q+1$, for
    all but a finite number of vertices $v\in VX$.  Then
 
    \itm{i} $ \chi_{av}(X)= \frac12(1-q)$ and
    \begin{align*}
	Z_{X,\cg}(u) &= (1-u^2)^{(1-q)/2}
	\Det_{\cg}((1+qu^2)I-uA)^{-1} \\
	&= (1-u^2)^{(1-q)/2} \Det_{\cg}((1+qu^2)I-(q+1)uP)^{-1},\
	\text{ for } |u|<\frac{1}{q},
    \end{align*}
  
    \itm{ii} by means of the determinant formula in $(i)$, $Z_{X,\cg}$
    can be extended to a function holomorphic (and without zeros) at
    least in the open set $\O$ defined in Lemma \ref{prop:holomorphy},
     
    \itm{iii} for $|u| < \frac1q$,
    $$
    \log Z_{X,\cg}(u) = \frac{1-q}{2} \log(1-u^2) +
    \sum_{n=1}^{\infty} \frac{1}{n} \sum_{k=0}^{n} \begin{pmatrix}
    n\\k \end{pmatrix} (q+1)^k (-q)^{n-k} u^{2n-k} Tr_\cg(P^k).
    $$
\end{Prop}
\begin{proof}
    Let us observe that, in the case of essentially $(q+1)$-regular
    graphs, we have $Tr_{\cg}(Q-I)= q-1$ so that, by Lemma
    \ref{lem:EuPoi}, the first part of $(i)$ follows and
    $$
    \frac{1}{Z_{X,\cg}(u)} = (1-u^{2})^{(q-1)/2}
    \Det_{\cg}((1+qu^{2})I-Au), \text{ for } |u|<\frac{1}{q}.
    $$
     
    $(i)$ Let $k\geq 0$ be the number of exceptional vertices, and
    assume that the vertices of $X$ have been ordered so that the
    first $k$ of them are the exceptional ones, then follow the
    vertices adjacent to the exceptional ones, and then all the
    others.  Then $A-(q+1)P = \d$, where only the first $k$ rows and
    the columns from $k+1$ to $k(q+1)$ of the matrix $\d$ can be
    nonzero.  Moreover, $Q=qI+\d'$, where $\d'$ is a diagonal matrix
    whose only possible nonzero entries are the first $k$ elements on
    the diagonal.  Therefore, $A-uQ = ((q+1)P-qu) + (\d-u\d')$ and
    \begin{align*}
	Tr_\cg((A-uQ)^n) & = \sum_{k=0}^{n} \begin{pmatrix} n\\k
	\end{pmatrix} Tr_\cg (((q+1)P-qu)^k(\d-u\d')^{n-k})\\
	& = Tr_\cg(((q+1)P-qu)^n) \\
	&\qquad + \sum_{k=0}^{n-1} \begin{pmatrix} n\\k \end{pmatrix}
	Tr_\cg (((q+1)P-qu)^k(\d-u\d')^{n-k}) \\
	& = Tr_\cg(((q+1)P-qu)^n),
    \end{align*}
    because 
    $$
    |Tr_\cg (((q+1)P-qu)^k(\d-u\d')^{n-k})| \leq \| (q+1)P-qu\|^k
    \|\d-u\d'\|^{n-k-1} Tr_\cg(|\d-u\d'|)
    $$
    and $Tr_\cg(|\d-u\d'|)=0$; indeed, only the first $k(q+1) \times
    k(q+1)$ block of the matrix $|\d-u\d'|$ can be nonzero.  The
    result follows.

    $(ii)$ From Lemma \ref{prop:holomorphy}, it follows that the
    factor $\Det_{\cg}((1+qu^2)I-uA)^{-1}$ is holomorphic in $\O$. 
    Then we define $(1-u^2)^{(1-q)/2}$ by means of its power series in
    the open disc $|u|<\frac{1}{q}$ and by holomorphic extension in
    $$
    \O':=\br^2 \setminus \left(\set{(x,y)\in\br^2: x^2+y^2=
    \frac{1}{q},\ y\leq 0} \cup \set{(x,0)\in\br^2: \frac{1}{q}\leq
    |x|\leq 1}\right).
    $$ 
    
    $(iii)$ This follows from $(ii)$.
\end{proof}

 \begin{rem}
     Observe that Proposition \ref{prop:essentiallyRegular} $(ii)$
     shows that, in the essentially regular case, the domain of
     holomorphicity $\O$ of $Z_{X,\cg}$ is larger than the open disc
     $\{u\in\bc: |u|<\frac{1}{d-1}\}$, which is the region where the
     Euler product is known to converge, as was proved in
     Theorem \ref{lem:power.series}.
 \end{rem}

 \begin{Thm} [Functional equations] \label{prop:functEqs}
     Let $X$ be essentially $(q+1)$-regular.  Then, for $u\in\O$,
     where $\O$ is the open set in Lemma
     \ref{prop:holomorphy}, we have
     
     \itm{i} $\La_{X}(u) := (1-u^{2})^{q/2} (1-q^{2}u^{2})^{1/2}
     Z_{X,\cg}(u) =-\La_{X}\Bigl(\frac{1}{qu}\Bigr)$,
     
     \itm{ii} $\xi_{X}(u) := (1+u)^{(q-1)/2} (1-u)^{(q+1)/2} (1-qu) 
     Z_{X,\cg}(u) = \xi_{X}\Bigl(\frac{1}{qu}\Bigr)$, 
     
     \itm{iii} $\Xi_{X}(u) := (1-u^{2})^{(q-1)/2} (1+qu^{2}) 
     Z_{X,\cg}(u) = \Xi_{X}\Bigl(\frac{1}{qu}\Bigr)$.
 \end{Thm}
 \begin{proof}
     $(i)$ 
     \begin{align*}
	 \La_{X}(u) & = (1-u^{2})^{1/2} (1-q^{2}u^{2})^{1/2}
	 \Det_{\cg}((1+qu^{2})I-Au)^{-1} \\
	 &= u\Bigl(\frac{q^{2}}{q^{2}u^{2}}-1 \Bigr)^{1/2} qu
	 \Bigl(\frac{1}{q^{2}u^{2}}-1 \Bigr)^{1/2}
	 \frac{1}{qu^{2}}\Det_{\cg}\Bigl(
	 (1+\frac{q}{q^{2}u^{2}})I-A\frac{1}{ qu} \Bigr)^{-1}\\
	 &= -\La_{X}\Bigl(\frac{1}{qu}\Bigr).
     \end{align*}
     $(ii)$
     \begin{align*}
	 \xi_{X}(u) & = (1-u) (1-qu) \Det_{\cg}((1+qu^{2})I-Au)^{-1}\\
	 & = u \Bigl( \frac{q}{qu} -1 \Bigr) qu \Bigl( \frac{1}{qu} -1
	 \Bigr) \frac{1}{qu^{2}}\Det_{\cg}\Bigl(
	 (1+\frac{q}{q^{2}u^{2}})I-A\frac{1}{ qu} \Bigr)^{-1}
	 = \xi_{X}\Bigl(\frac{1}{qu}\Bigr).
     \end{align*}
     $(iii)$
     \begin{align*}
	 \Xi_{X}(u) & = (1+qu^{2}) \Det_{\cg}((1+qu^{2})I-Au)^{-1} \\
	 &= qu^{2} \Bigl( \frac{q}{q^{2}u^{2}} +1 \Bigr)
	 \frac{1}{qu^{2}}\Det_{\cg}\Bigl(
	 (1+\frac{q}{q^{2}u^{2}})I-A\frac{1}{ qu} \Bigr)^{-1}
	 = \Xi_{X}\Bigl(\frac{1}{qu}\Bigr).
     \end{align*}
 \end{proof}
 
\begin{exmp}\label{exmp:chi(gasket)}
    Let $X$ be the Gasket graph given in figure
    \ref{fig:GasketVicsek}a.  Then it is clear that, except for the
    marked point in figure \ref{fig:GasketVicsek}a, all the vertices
    have degree $4$.  Hence, $X$ is essentially $4$-regular, and
    Proposition \ref{prop:essentiallyRegular} and Theorem
    \ref{prop:functEqs} apply, with $q=3$.  Entirely analogous
    examples are provided by the higher-dimensional ``gaskets'' based
    on the higher-dimensional simplexes.  On the other hand, it can be
    easily checked that the other examples of self-similar graphs
    given in figures \ref{fig:GasketVicsek}b,
    \ref{fig:LindstromCarpet}a and \ref{fig:LindstromCarpet}b are not
    essentially regular.
\end{exmp}

\section{Approximation by finite graphs}\label{sec:approx}

 In this final section, we obtain an approximation result which shows 
 that our definition of the zeta function is a natural one.
 
 \begin{Lemma}\label{lem:normEstimate}
     Let $X$ be a self-similar graph, as in Section
     \ref{sec:selfSimilar}.  Let $A$ be the adjacency operator as in
     (\ref{eq:adjacency}) and $Q=D-I$ where $D$ is as in
     (\ref{eq:degree}).  Let $f(u) := Au-Qu^2$, for $u\in\bc$.  Then
     $\|f(u)\|<\frac12$, for $|u|<\frac{1}{d+\sqrt{d^2+2(d-1)}}$.
 \end{Lemma}
 \begin{proof}
     Since $\displaystyle
     \|f(u)\|  \leq |u| \|A\| + |u|^2\|Q\| \leq d|u|+(d-1)|u|^2$
     for any $u\in\bc$, the thesis follows.
 \end{proof}

 \begin{Thm}[Approximation by finite graphs]
     Let $X$ be a self-similar graph, as in Section 
     \ref{sec:selfSimilar}. Then
     $$
     Z_{X,\cg}(u) = \lim_{n\to\infty} Z_{K_{n}}(u)^{\frac{1}{|K_{n}|}},     
     $$
     uniformly on compact subsets of $\set{ u\in\bc:
     |u| <\frac{1}{d+\sqrt{d^2+2(d-1)}} }$.
 \end{Thm}
 \begin{proof}
     Let  $f(u):= Au-Qu^{2}$ and $E_{n}:= P(K_{n})$. Since the $Q$-matrix $Q_n$ for the sub-graph $K_n$ is different form $E_nQE_n$, we set $ \delta_n=E_nQE_n -Q_n$, which is a positive diagonal matrix. Then, 
     with $Tr$ denoting the usual trace on bounded operators, we have
     $$
     \log Z_{K_{n}}(u) = -\frac12 Tr(E_{n}(Q-I)E_{n}) \log(1-u^{2}) -
     Tr \log(E_{n}(I-f(u))E_{n})+o(|K_n|), \  n\to\infty.
     $$
     In fact, 
     \begin{align*}
	\log Z_{K_n}(u) = & -\frac12 Tr(E_n(Q-I)E_n) \log(1-u^2) + \frac12 Tr(\delta_n) \log(1-u^2) \\
	& - Tr \log(E_n(I- f(u)-\delta_n u^2)E_n ).
\end{align*}
     Now $Tr \log(E_n(I- f(u)-\delta_n u^2)E_n ) = -\sum_{k=1}^\infty \frac1k Tr( (E_nf(u)E_n +\delta_n u^2)^k)$, and $Tr( (E_nf(u)E_n +\delta_n u^2)^k) = Tr( (E_n f(u) E_n)^k) + a_{kn}$, where $a_{kn}$ is the sum of $2^k-1$ terms, each being the trace of a product of $k$ operators, at least one of them being $\delta_n$, so that $|a_{kn}| \leq \sum_{j=0}^{k-1} \binom{k}{j} \| f(u) \|^j (d-1)^{k-j} |u|^{2(k-j)} \eps_n |K_n| \leq ( \|f(u)\| +(d-1)|u|^2)^k \eps_n |K_n| \leq (d|u|+2(d-1)|u|^2)^k \eps_n |K_n|$. Therefore, if $|u|< \frac{1}{d+\sqrt{d^2+2(d-1)}}$,  
     $$
     \frac1{|K_n|} \bigl| Tr \log(E_n(I- f(u)-\delta_n u^2)E_n ) - Tr \log(E_n(I- f(u))E_n ) \bigr| \leq \eps_n \log \frac{1}{d|u|+2(d-1)|u|^2} \to 0.
     $$
     Since $\frac{1}{|K_n|} Tr(\delta_n) \to 0$, we proved our statement.
     
     \noindent Moreover,
     $$
     Tr \log(E_{n}(I-f(u))E_{n}) = -\sum_{k=1}^{\infty} 
     \frac{1}{k} Tr\bigl( (E_{n}f(u)E_{n})^{k}\bigr).
     $$
     Observe that, for $k\geq 2$,
     \begin{align*}
	 Tr \bigl( E_{n}f(u)^{k}E_{n} \bigr) &= Tr\bigl( E_{n} (
	 f(u)(E_{n}+E_{n}^{\perp}) )^{k} E_{n} \bigr)\\
	 &= Tr\bigl( (E_{n}f(u)E_{n})^{k} \bigr) + \sum_{ \substack{
	 \s\in\{-1,1\}^{k-1}\\ \s\neq \{1,1,\ldots,1\} } } Tr \bigl(
	 E_{n} \prod_{j=1}^{k-1} [f(u)E_{n}^{\s_{j}}]f(u)E_{n} \bigr),
     \end{align*} 
     where $E_n^{-1}$ stands for $E_n^\perp$, and
     \begin{align*}
	 | Tr \bigl( E_{n} \prod_{j=1}^{k-1} [f(u)E_{n}^{\s_{j}}]
	 f(u)E_{n} \bigr) | & = | Tr \bigl( ...E_{n} f(u)E_{n}
	 ^\perp...  \bigr) | \\
	 & \leq \| f(u) \|^{k-1} Tr ( |E_{n} f(u)E_{n}^\perp|).
     \end{align*}	 
     Moreover, with $\O_{n}:= B_{1}(VK_{n})\setminus VK_n \subset 
     B_{1}(\cf_{\cg} K_{n})$, we have
     \begin{align*}
	 Tr ( |E_{n} f(u)E_{n}^\perp|) & = Tr (|P(K_{n} )
	 f(u)P(\O_{n})| )\\
	 & \leq \|f(u)\| Tr(P(\O_{n})) \\
	 & = \|f(u)\| |\O_{n}| \\
	 & \leq \|f(u)\| (d+1) \eps_n |K_n|.
     \end{align*}	 
     Hence,   we obtain
     $$
     | Tr \bigl( E_{n}f(u)^{k}E_{n} \bigr) - Tr\bigl(
     (E_{n}f(u)E_{n})^{k} \bigr) | \leq (2^{k-1}-1) \| f(u) \|^{k} (d+1)
     \eps_n |K_n|,
     $$
     so that 
     \begin{align*}
	 \biggl| Tr&\log(E_{n}(I-f(u))E_{n}) -
	 Tr(E_{n} \log(I-f(u)) E_{n}) \biggr| \\
	 &= \biggl| \sum_{k=1}^{\infty} \frac{1}{k} Tr\bigl(
	 (E_{n}f(u)E_{n})^{k} \bigr) - \sum_{k=1}^{\infty}
	 \frac{1}{k} Tr \bigl( E_{n}f(u)^{k}E_{n} \bigr)\biggr| \\
	 & \leq \biggl( \sum_{k=1}^{\infty}
	 \frac{2^{k-1}\|f(u)\|^{k}}{k}\biggr)(d+1) \eps_n |K_n| \\
	 & \leq C (d+1) \eps_n |K_n|,
     \end{align*}
     where the series converges for $|u|<
     \frac{1}{d+\sqrt{d^2+2(d-1)}}$, by Lemma \ref{lem:normEstimate}. 
     Therefore,
     $$
     \biggl| \frac{Tr\log(E_{n}(I-f(u))E_{n})}{|K_{n}|} -
     \frac{Tr(E_{n} \log(I-f(u)) E_{n})}{|K_{n}|} \biggr| \to 0, \
     \text{ as } n\to\infty,
     $$
     and 
     \begin{equation*}
	\begin{split}
	 \lim_{n\to\infty}& \frac{\log Z_{K_{n}}(u)}{|K_{n}|} = \\
	 &= -\frac12 
	 \lim_{n\to\infty} \frac{Tr(E_{n}(Q-I)E_{n}) }{|K_{n}|}
	 \log(1-u^{2})  - \lim_{n\to\infty} \frac{Tr(E_{n} \log(I-f(u)) 
	 E_{n})}{|K_{n}|} \\
	 &= -\frac12 Tr_{\cg}(Q-I) \log(1-u^{2}) -
	 Tr_{\cg}(\log(I-f(u)))  = \log Z_{X,\cg}(u),
	 \end{split}
     \end{equation*}
     from which the claim follows. 
 \end{proof}

 \begin{rem}
     Observe that $\frac{1}{2\a}<\frac{1}{d+\sqrt{d^2+2(d-1)}}
     <\frac{1}{\a}$.
 \end{rem} 
 
 \medskip\par\noindent{\it Acknowledgements.}
     The second and third named authors would like to thank
     respectively the University of California, Riverside, and the
     University of Roma ``Tor Vergata'' for their hospitality at
     various stages of the preparation of this paper.
 

\end{document}